\documentstyle[twoside,psfig,amssymb]{article}

\oddsidemargin=\evensidemargin 
\catcode`\@=11 \long\def\@makefntext#1{ \protect\noindent \hbox to
3.2pt {\hskip-.9pt
$^{{\eightrm\@thefnmark}}$\hfil}#1\hfill}       

\def\ps@myheadings{\let\@mkboth\@gobbletwo      
\def\@oddhead{\hbox{}
\rightmark\hfil\eightrm\thepage}
\def\@oddfoot{}\def\@evenhead{\eightrm\thepage\hfil
\leftmark\hbox{}}\def\@evenfoot{}
\def\sectionmark##1{}\def\subsectionmark##1{}}

\catcode`@=11                        
\def\ps@plain{\let\@mkboth\@gobbletwo
     \def\@oddhead{}\def\@oddfoot{\eightrm\hfil\thepage
     \hfil}\def\@evenhead{}\let\@evenfoot\@oddfoot}

\newcounter{sectionc}\newcounter{subsectionc}\newcounter{subsubsectionc}
\renewcommand{\section}[1] {\vspace{12pt}\addtocounter{sectionc}{1}
\setcounter{subsectionc}{0}\setcounter{subsubsectionc}{0}\noindent
    {\tenbf\thesectionc. #1}\par\vspace{5pt}}
\renewcommand{\subsection}[1] {\vspace{12pt}\addtocounter{subsectionc}{1}
    \setcounter{subsubsectionc}{0}\noindent
    {\bf\thesectionc.\thesubsectionc.
    {\kern1pt \bfit #1}}\par\vspace{5pt}}
\renewcommand{\subsubsection}[1] {\vspace{12pt}
    \addtocounter{subsubsectionc}{1}
    \noindent
    {\tenrm\thesectionc.\thesubsectionc.\thesubsubsectionc. {\kern1pt
    \it #1}}\par\vspace{5pt}}

\newcounter{appendixc}
\newcounter{subappendixc}[appendixc]
\newcounter{subsubappendixc}[subappendixc]

\renewcommand{\appendix}[1] {\vspace{12pt}  
    \refstepcounter{appendixc}      
    \setcounter{figure}{0}
    \setcounter{table}{0}
    \setcounter{lemma}{0}
    \setcounter{theorem}{0}
    \setcounter{corollary}{0}
    \setcounter{definition}{0}
    \setcounter{equation}{0}
    \renewcommand{\thefigure}{\Alph{appendixc}.\arabic{figure}}
    \renewcommand{\thetable}{\Alph{appendixc}.\arabic{table}}
    \renewcommand{\theappendixc}{\Alph{appendixc}}
    \renewcommand{\thelemma}{\Alph{appendixc}.\arabic{lemma}}
    \renewcommand{\thetheorem}{\Alph{appendixc}.\arabic{theorem}}
    \renewcommand{\thedefinition}{\Alph{appendixc}.\arabic{definition}}
    \renewcommand{\thecorollary}{\Alph{appendixc}.\arabic{corollary}}
    \renewcommand{\theequation}{\Alph{appendixc}.\arabic{equation}}
    \noindent{\tenbf Appendix \theappendixc #1}\par\vspace{5pt}}

\newcommand{\smalllineskip}{\baselineskip=10pt}

\newcommand{\copyrightheading}[1]
    {\vspace*{-2.5cm}\smalllineskip{\flushleft
    {\footnotesize }\\
    {\footnotesize \copyright\kern2pt }\\
         }}

\def\keywords#1{{
    \centering{\begin{minipage}{4.5in}\footnotesize\baselineskip=10pt
    {\footnotesize\it Keywords}\/: #1
    \end{minipage}}\par}}

\newcounter{itemlistc}
\newcounter{romanlistc}
\newcounter{alphlistc}
\newcounter{arabiclistc}

\newcommand{\fcaption}[1]{
        \refstepcounter{figure}
        \setbox\@tempboxa = \hbox{\footnotesize Fig.~\thefigure. #1}
        \ifdim \wd\@tempboxa > 5in
           {\begin{center}
        \parbox{5in}{\footnotesize\smalllineskip Fig.~\thefigure. #1}
            \end{center}}
        \else
             {\begin{center}
             {\footnotesize Fig.~\thefigure. #1}
              \end{center}}
        \fi}

\newcommand{\tcaption}[1]{
        \refstepcounter{table}
        \setbox\@tempboxa = \hbox{\footnotesize Table~\thetable. #1}
        \ifdim \wd\@tempboxa > 5in
           {\begin{center}
        \parbox{5in}{\footnotesize\smalllineskip Table~\thetable. #1}
            \end{center}}
        \else
             {\begin{center}
             {\footnotesize Table~\thetable. #1}
              \end{center}}
        \fi}

\def\pmb#1{\setbox0=\hbox{#1}
    \kern-.025em\copy0\kern-\wd0
    \kern.05em\copy0\kern-\wd0
    \kern-.025em\raise.0433em\box0}

\def\fnt#1#2{\footnotetext{\kern-.3em
    {$^{\mbox{\scriptsize #1}}$}{#2}}}

\def\fpage#1{\begingroup
\voffset=.3in
\thispagestyle{empty}\begin{table}[b]\centerline{\footnotesize #1}
    \end{table}\endgroup}

\def\runninghead#1#2{\pagestyle{myheadings}
\markboth{{\protect\footnotesize\it{\quad #1}}\hfill}
{\hfill{\protect\footnotesize\it{#2\quad}}}}

\font\tenrm=cmr10  \font\tenbf=cmbx10
\font\bfit=cmbxti10 at 10pt \font\ninerm=cmr9 
 \font\eightrm=cmr8

\newtheorem{theorem}{Theorem}   

\newtheorem{lemma}{Lemma}

\newtheorem{definition}{Definition}

\def\@begintheorem#1#2{\trivlist    
    \item[\hskip\labelsep{\bf #1\ #2.}]}
\def\@opargbegintheorem#1#2#3{\trivlist
    \item[\hskip\labelsep{\bf #1\ #2\ (#3).}]}

        {\setcounter{itemlistc}{0}      
     \begin{list}{$\bullet$}        
    {\usecounter{itemlistc}         
     \leftmargin10pt           
     \setlength{\parsep}{0pt}
     \setlength{\itemsep}{0pt}     
    }}{\end{list}}

    {\setcounter{romanlistc}{0}     
     \begin{list}{$($\roman{romanlistc}$)$} 
    {\usecounter{romanlistc}        
     \leftmargin18pt 
     \setlength{\parsep}{0pt}
     \setlength{\itemsep}{0pt}  
     \settowidth{\labelwidth}{#1}
    }}{\end{list}}

    {\setcounter{enumii}{0}         
     \begin{list}{$($\alph{enumii}$)$}  
    {\usecounter{enumii}            
     \leftmargin18pt        
     \setlength{\parsep}{0pt}
     \setlength{\itemsep}{0pt}  
     \settowidth{\labelwidth}{#1}
    }}{\end{list}}

\def\qed{\hbox{${\vcenter{\vbox{            
   \hrule height 0.4pt\hbox{\vrule width 0.4pt height 6pt
   \kern5pt\vrule width 0.4pt}\hrule height 0.4pt}}}$}}

\def\theequation{\thesectionc.\arabic{equation}}  
\begin{document}

\runninghead{Adequacy of link families } {Adequacy  of link
families}

\setcounter{page}{1}

\markboth{Slavik Jablan}{}





\fpage{1} \centerline{\bf ADEQUACY OF LINK FAMILIES}
\bigskip

\centerline{\footnotesize SLAVIK JABLAN}
\medskip
\centerline{\footnotesize\it The Mathematical Institute, Knez
Mihailova 36,}\centerline{\footnotesize\it P.O.Box 367, 11001
Belgrade,}\centerline{\footnotesize\it Serbia}
\centerline{\footnotesize\it jablans@mi.sanu.ac.yu}

\bigskip

\bigskip
\bigskip

\begin{abstract}
Using computer calculations and working with representatives of
pretzel tangles we established general adequacy criteria for
different classes of knots and links. Based on adequate graphs
obtained from all Kauffman states of an alternating link we defined
a new numerical invariant: adequacy number, and computed adequacy
polynomial which is the invariant of alternating link families.
Adequacy polynomial distinguishes (up to mutation) all families of
alternating knots and links whose generating link has at most $n=12$
crossings.
\end{abstract}


\section{Introduction}

First we give a brief overview of the properties of adequate,
semi-adequate and inadequate link diagrams and their corresponding
links. In this paper, we will consider only prime links.

Let $D$ be a diagram of an unoriented link $L$ framed in a 3-ball
$B^3$. A Kauffman state of a diagram $D$ is a function from the set
of crossings of $D$ to the set of signs $\{+1,-1\}$. Graphical
interpretation is smoothing each crossing of $D$ by introducing
markers according to the convention illustrated in Fig. 1. A {\it
state diagram} $D_s$ is a system of circles obtained by smoothing
all crossings of $D$ [PrAs]. The set of circles in $D_s$, which are
called {\it state circles}, is denoted by $C(D)$. Points of the
state circles corresponding to a smoothed crossing are called {\it
touch-points}. The number of touch-points belonging to a state
circle $c\in C(D)$ is called the {\it length} of $c$.

\begin{figure}[th]
\centerline{\psfig{file=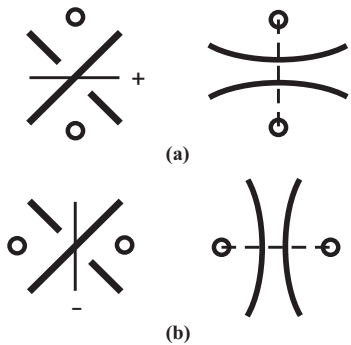,width=2.00in}}
\vspace*{8pt} \caption{(a) $-$marker; (b) $+$marker. The broken
lines represent the edges of the associated graph $G_s$ connecting
state circles (represented by dots). \label{f1.1}}
\end{figure}

Kauffman states $s_+$ and $s_-$ with all $+$ or all $-$ signs are
called {\it special states}, and their corresponding state
diagrams $D_{s_+}$ and $D_{s_-}$ are called {\it special
diagrams}. All other Kauffman states with both $+$ or $-$ signs
are called {\it mixed states}, and to them correspond {\it mixed
state diagrams}.

\begin{definition}
A diagram $D$ is $s$-adequate if two arcs at every touch-point of
$D_s$ belong to different state circles. In particular, a diagram
$D$ is $+$adequate or $-$adequate if it is $s_+$ or $s_-$
adequate, respectivelly.  If a diagram is neither $+$adequate nor
$-$adequate it is called {\it inadequate}. If a diagram is both
$+$adequate and $-$adequate, it is called {\it adequate}, and if
it is only $+$adequate or $-$adequate, it is called {\it
semi-adequate} [LiThi,Li].
\end{definition}

To every state diagram $D_s$ we associate the graph $G_{s}$, whose
vertices are state circles of $D_s$ and edges are lines connecting
state circles via smoothed crossings in $D$. Now we can restate
Definition 1 in terms of $G_s$: $D$ is $s$-adequate if $G_s$ is
loopless. A state graph $G_s$ is called adequate if $D_s$ is
$s$-adequate.

\begin{definition}
A link is {\it adequate} if it has an adequate ($+$adequate and
$-$adequate) diagram. A link is {\it semi-adequate} if it has a
$+$ or $-$adequate diagram. A link is {\it inadequate} if it is
neither $+$ or $-$adequate [LiThi,Li].
\end{definition}

The mirror image of a diagram transforms the $+$adequacy into
$-$adequacy and {\it vice versa}.

\begin{definition}
A link that has one $+$adequate diagram and another diagram that
is $-$adequate is called {\it weakly adequate}.
\end{definition}

\noindent For example, knot $11n_{146}$ $9^*.-2:.-2$ has
$-$adequate 11-crossing diagram and $+$adequate 12-crossing
diagram $6^*-2.2.-2.2.2\,0.-2\,0$. Another such example is Perko's
knot $10_{161}$ $3:-2\,0:-2\,0$ (Fig. 4) [Stoi].

A crossing in a link diagram for which there exists a circle in
the projection plane intersecting the diagram transversely at that
crossing, but not intersecting the diagram at any other point is
called {\it nugatory} crossing. A link diagram is called {\it
reduced} if it has no nugatory crossings. The following theorem
holds for reduced alternating link diagrams:

\begin{theorem}
A reduced alternating diagram is adequate [LiThi,Li,Cro].
\end{theorem}

\noindent Hence, all alternating links are adequate.

\begin{theorem}
An adequate diagram has minimal crossing number [LiThi,Li,Cro].
\end{theorem}

\noindent This theorem can be used to prove minimality of some
non-alternating link diagrams.

\begin{theorem}
Every unlink diagram is inadequate. Semi-adequate link diagrams
are non-trivial [Thist].
\end{theorem}

A non-minimal diagram of an adequate link can be semi-adequate or
inadequate. For example, non-minimal diagram $3\,2\,4\,-2\,2$ of
the alternating knot $3\,3\,2\,3$ is semi-adequate, and
non-minimal diagram $3\,3\,4\,-1\,2$ of the alternating knot
$3\,3\,2$ is inadequate.

A non-minimal diagram of a semi-adequate link also can be
semi-adequate or inadequate. For example, non-minimal diagram
$3,3,2,2\,-3$ and minimal diagram of the same knot $3,3,2,-2\,-2$
are both semi-adequate; minimal diagram of the knot $2\,1,3,-2$ is
semi-adequate, and it's non-minimal diagram $2\,1,3,2-$ is
inadequate.

\begin{theorem}
Two adequate diagrams of a link have the same crossing number and
the same writhe [Cro].
\end{theorem}

\begin{definition}
An alternating diagram of a marked 2-tangle $t$ is called {\it
strongly alternating} if the both its closures, numerator closure
$N(t)$ and denominator closure $D(t)$, are irreducible
[LiThi,Li,Cro].
\end{definition}

\begin{theorem}
The non-alternating sum of two strongly alternating tangles is
adequate [LiThi,Li,Cro].
\end{theorem}

\noindent This theorem can be very efficiently used to prove that
certain types of link diagrams are adequate. For example, all
semi-alternating diagrams are adequate [LiThi,Li]

According to Theorem 2, minimal diagrams can be used to determine if
a link is adequate, but do not provide necessary and sufficient
conditions to distinguish semi-adequate links from inadequate ones.

\begin{theorem}
A link is inadequate if both coefficients of the terms of highest
and lowest degree of its Jones polynomial are different from $\pm
1$.
\end{theorem}

The proof of this theorem for knots follows directly from the
results of W.B.R.~Lickorish and M.~Thistlethwaite, and it also holds
for links, due to work of J.~Przytycki [LiThi,Pr].

\section{Adequate links with at most 12 crossings}

Using {\it Knotscape} tables of knots given in Dowker-Thistlethwaite
notation, A.~Stoi- menow detected all non-alternating adequate knots
up to $n=16$ crossings. In this paper we consider adequacy of
non-alternating links and their families (classes) given in Conway
notation.

Adequate non-alternating links with $n\le 10$ crossings are given
in the following table:

\small

\bigskip

\begin{tabular}{|c|c|c|} \hline
 $n=8$ & & \\ \hline
  $2,2,-2,-2$ & $(2,2) -(2,2)$ &  \\ \hline
 2 Links & & \\ \hline \hline

$n=9$ & & \\ \hline

$3,2,-2,-2$ & $(3,2)\,-(2,2)$ & $(2\,1,2)\,-(2,2)$ \\ \hline

$.-(2,2)$ & & \\ \hline

4 Links & & \\ \hline \hline

 $n=10$ & & \\ \hline

$(3,2)\,-(3,2)$  & $(3,2)\,-(2\,1,2)$ & $(2\,1,2)\,-(2\,1,2)$ \\
\hline

 3 Knots & & \\ \hline

 $3,2\,1,-2,-2$  & $3,3,-2,-2$ & $3,-2,2\,1,-2$ \\ \hline

 $3,-2,3,-2$  & $4,2,-2,-2$ & $2,2,2,-2,-2$ \\ \hline

 $2\,2,2,-2,-2$  & $(4,2)\,-(2,2)$ &  $(3,2\,1)\,-(2,2)$ \\ \hline

 $(3\,1,2)\,-(2,2)$  & $(2\,1,2\,1)\,-(2,2)$ &  $(3,3)\,-(2,2)$ \\ \hline

 $(2\,1\,1,2)\,-(2,2)$  & $(2,-2,-2)\,(2,2)$ & $(2\,2,2)\,-(2,2)$ \\ \hline

$(2,2,2)\,-(2,2)$ & $(2,2),2,-(2,2)$ & $.-(2,2).2$ \\ \hline

$.-(2,2).2\,0$ & $.-(2,2):2\,0$ & $.-(2,2):2$ \\ \hline

$103^*-1.-1.-1.-1::.-1$ &   &  \\ \hline

 22 Links & & \\ \hline

\end{tabular}

\normalsize

\bigskip

\noindent All of them, except polyhedral ones, satisfy Theorem 5
or are obtained from the pretzel links which satisfy this theorem
by permuting their rational tangles.

Theorem 6 gives sufficient but not necessary conditions for
recognizing inadequate links. For example, the first and last
coefficient of Jones polynomial of the knot $11n_{95}$ =
$2\,0.-2\,1.-2\,0.2$ are different from $\pm 1$, so it is
inadequate [Cro]. However, since this theorem doe's not give
necessary conditions for a link to be inadequate, the main problem
remains detection of inadequate links.

For knots with at most $n\le 12$ crossings every minimal diagram
of a semi-adequate knot is semi-adequate. Unfortunately, this is
not true for knots with $n\ge 13$ crossings: the first example of
a semi-adequate knot with a minimal inadequate diagram (Fig. 2) is
the knot $13n_{4084}$ $10^{**}.-1.-1.-1:.-2.2.-2$ with the minimal
Dowker-Thistlethwaite code $$\{\{13\}, \{6, -10, 12, 24, 20, -18,
-26, -22, -4, 2, -16, 8, -14\}\}.$$ \noindent Except this
inadequate diagram of writhe 9, it has another semi-adequate
minimal diagram $11^{**}.-2::-2\, 0:-1.-1.-1$ of writhe 7, with
the Dowker-Thistlethwaite code
$$\{\{13\}, \{6, 12, -16, 23, 2, 17, 21, 26, 11, -4, -25, 7,
13\}\}$$ [KidSto, Stoi2]. For $n=15$ appear first semi-adequate
knots without a minimal semi-adequate diagram. For example, knot
$15n_{164563}$ has only minimal diagram
$10^{**}-1.-2\,0.2\,0::.2\,0.2\,0.-2\,0$, and it is inadequate
(Fig. 3). However, it has 16-crossing diagram
$11^*2\,0.-1.-2.-1.3\, 0.-1.2\,0::-1$ which is semi-adequate
[Stoi3]. This example can be generalized to the family of knot
diagrams $10^{**}-1.-2\,0.(2k)\,0::.2\,0.2\,0.-2\,0$ and
$11^*(2k)\,0.-1.-2.-1.3\, 0.-1.2\,0::-1$ ($k\ge 1$) with the same
properties, respectively.

\begin{figure}[th]
\centerline{\psfig{file=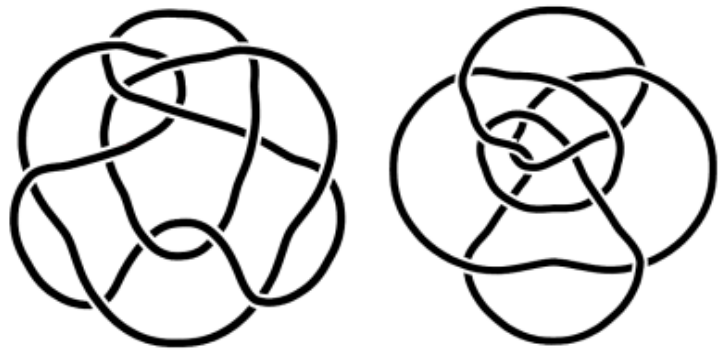,width=4.00in}} \vspace*{8pt}
\caption{Semi-adequate knot $13n_{4084}$ with a minimal inadequate
diagram $10^{**}.-1.-1.-1:.-2.2.-2$ and minimal semi-adequate
diagram $11^{**}.-2::-2\, 0:-1.-1.-1$ [KidSto,Stoi2].
\label{f1.1}}
\end{figure}

\begin{figure}[th]
\centerline{\psfig{file=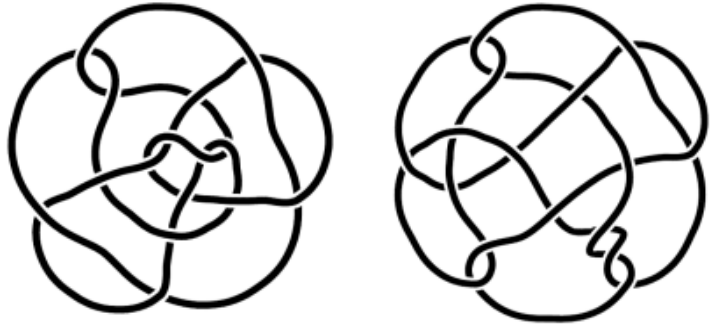,width=4.00in}} \vspace*{8pt}
\caption{Semi-adequate knot $15n_{164563}$ which has only minimal
diagram $10^{**}-1.-2\,0.2\,0::.2\,0.2\,0.-2\,0$ which is
inadequate and non-minimal 16-crossing diagram
$11^*2\,0.-1.-2.-1.3\, 0.-1.2\,0::-1$ which is semi-adequate
[Stoi3]. \label{f1.1}}
\end{figure}

For knots with at most $n\le 12$ crossings we checked adequacy
using all their minimal diagrams, but for all links and knots with
$n\ge 13$ crossings for each link or knot we used only one minimal
diagram.

The sign of adequacy is not necessarily the same for all minimal
diagrams of the same link, so we obtain weakly adequate links.

An example of a weakly adequate knot is Perko pair
$6^*3:-2\,0:-2\,0$ and $6^*-2\,-1.-1.2\,0.-1.2\,0.-1$ [Stoi]. This
example generalizes to one-parameter knot families called {\it
Perko families} [JaSaz]. Conway symbols $6^*(2k+1):-2\,0:-2\,0$
and $6^*-(2k)\,-1.-1.2\,0.-1.2\,0.-1$ represent two families of
minimal diagrams of the same weakly adequate knots with adequacy
of different signs and different writhe. For $k=1$ we obtain Perko
pair (Fig. 4), for $k=2$ two diagrams of the knot $12n_{850}$, for
$k=3$ two diagrams of the knot $14n_{26229}$, and for $k=4$ two
diagrams of the knot $16n_{965076}$ given in {\it Knotscape}
notation. The same holds for the minimal diagrams
$6^*2\,(2k):-2\,0:-2\,0$ and
$6^*-2\,-(2k-1)\,-1.-1.2\,0.-1.2\,0.-1$ of the knots $11n_{135}$,
$13n_{3546}$, and $15n_{114094}$ obtained for $k=1,2,3$,
respectively. Hence, for every $n\ge 10$ there exists at least one
weakly adequate knot which has two minimal diagrams with adequacy
of different signs and different writhe. Moreover, if $t$ is any
positive rational tangle ($t\neq 1$)\footnote{A rational tangle is
called positive if its Conway symbol contains only positive
numbers, and negative if it contains only negative numbers.},
minimal diagrams $6^*t\,(k+1):-2\,0:-2\,0$ and
$6^*(-t)\,(-k)\,(-1).-1.2\,0.-1.2\,0.-1$ of the same link have
adequacy of different signs and different writhe. Two minimal
diagrams of the knot obtained for $t=2\,2$ and $k=3$ are
illustrated in Fig. 5. In all these cases, the writhe of the
diagrams differs by 2, the first diagram is $+$adequate, and the
other $-$adequate. Since the class $6^*t\,(k+1):-2\,0:-2\,0$
contains links as well (e.g., for $t=2\,1$, $k=2$), this is the
first example of weakly adequate links.

\begin{figure}[th]
\centerline{\psfig{file=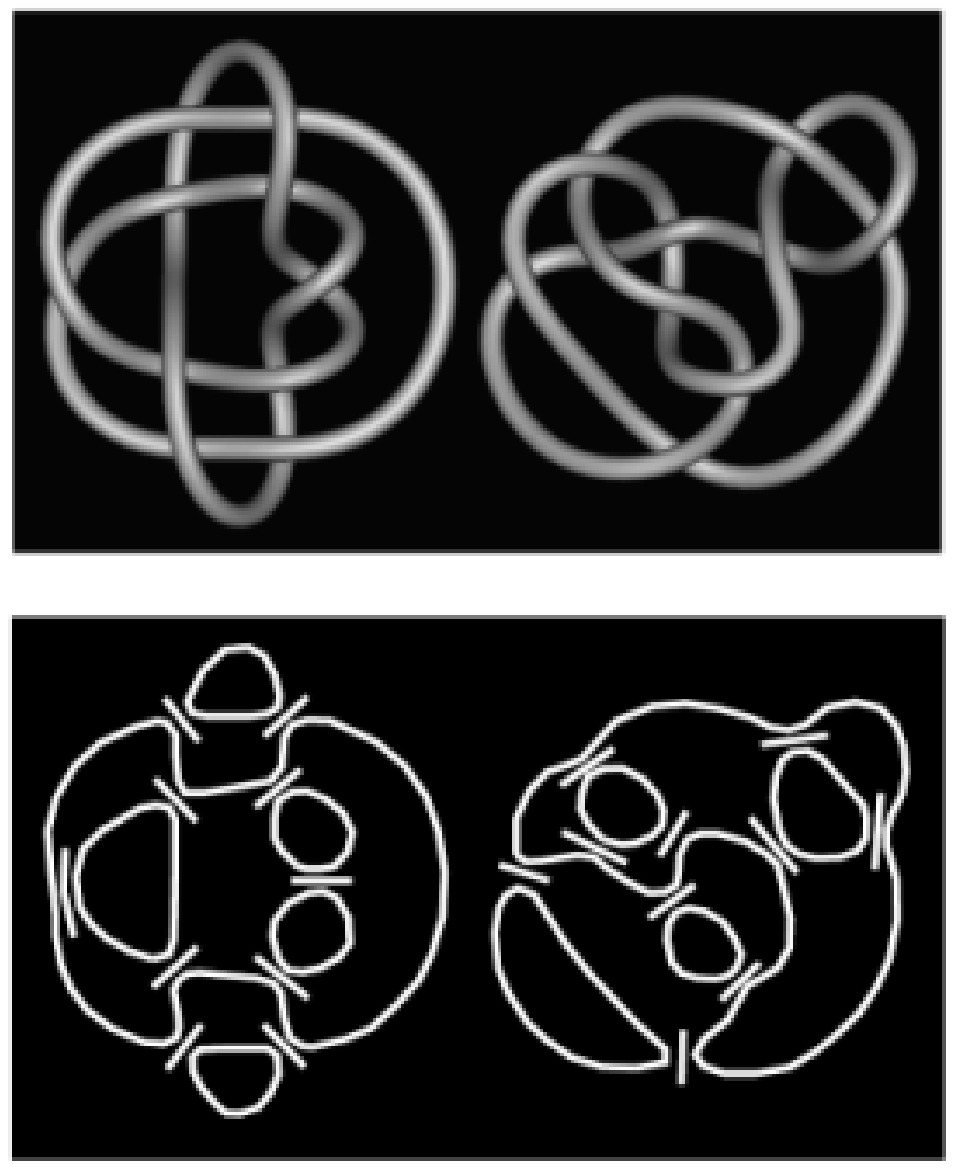,width=3.00in}}
\vspace*{8pt} \caption{Perko pair: semi-adequate knot with two
minimal diagrams $6^*3:-2\,0:-2\,0$ and
$6^*-2\,-1.-1.2\,0.-1.2\,0.-1$ with the adequacy of different
signs, where the first is $+$adequate, and the other
$-$adequate.\label{f1.1}}
\end{figure}

\begin{figure}[th]
\centerline{\psfig{file=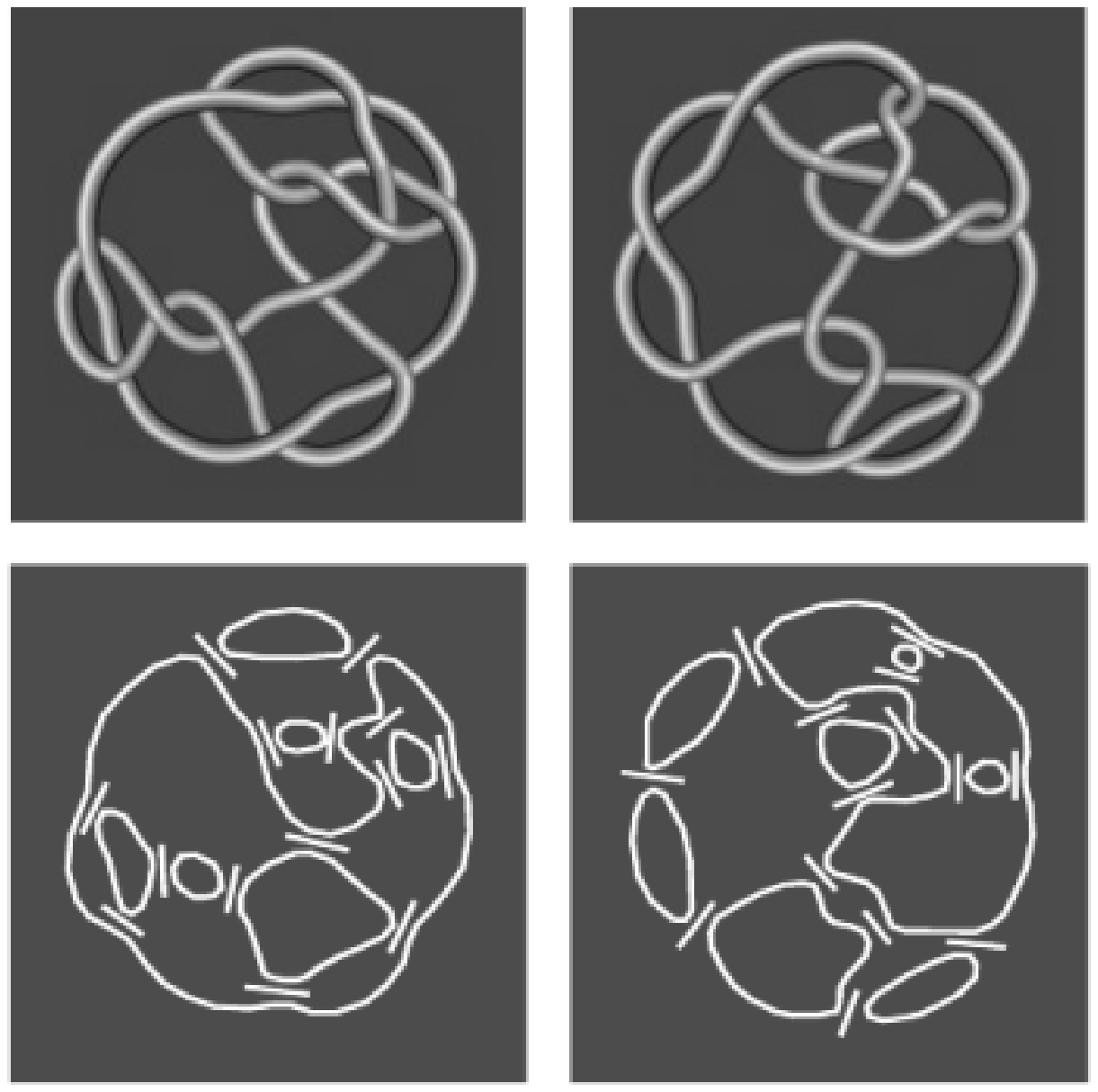,width=4.00in}} \vspace*{8pt}
\caption{Perko-type pair of knot diagrams: semi-adequate knot with
two minimal diagrams $6^*2\,2\,3:-2\,0:-2\,0$ and
$6^*-2\,-2\,-2\,-1.-1.2\,0.-1.2\,0.-1$ with the adequacy of
different signs, where the first is $+$adequate, and the other
$-$adequate. \label{f1.1}}
\end{figure}

At least for small number of crossings, most of non-alternating
links are semi-adequate, so adequate and inadequate links
represent a small portion of all non-alternating links. Hence, it
is of interest to tabulate adequate non-alternating links and
candidates for inadequate links and try to find some general
criteria for adequacy of certain classes of links. We checked
adequacy of all minimal diagrams of non-alternating knots and
links with at most $n=12$ crossings given in Conway notation.

Among 202 non-alternating links with at most 10 crossings there
are only 28 adequate links and 3 adequate knots. Links with
inadequate minimal diagrams are even more rare. Their list for
$n=10$ is given in the following table:

\bigskip

\small

\begin{tabular}{|c|c|c|} \hline
 $n=10$ & & \\ \hline
 $(2,2,-2)\,(2,-2)$ & $2.-2\,0.-2.2\,0$ & $103^*-1.-1::-1.-1$ \\ \hline
 3 Links & & \\ \hline

\end{tabular}

\normalsize

\bigskip

\noindent where links $2.-2\,0.-2.2\,0$ and $103^*-1.-1::-1.-1$
are inadequate according to Theorem 6, and nothing is known for
the link $(2,2,-2)\,(2,-2)$.

Particular links, families or classes of links which have a
minimal inadequate diagram will be refered to us as {\it
candidates for inadequate links} and in some cases Theorem 6 will
confirm that they indead are inadequate.

Candidates for inadequate knots occur for the first time among
11-crossing knots: knot $2\,0.-2\,1.-2\,0.2$ is inadequate
according to Theorem 6, but for the knot $2\,0.-3.-2\,0.2$ which
all minimal diagrams are inadequate, it is not possible to make
any conclusion, since both leading coefficients of its Jones
polynomial are equal to 1.

For $n=12$, among 19 knots with an inadequate minimal diagram, 11
knots given in the following table are inadequate according to
Theorem 6

\bigskip

\small

\begin{tabular}{|c|c|c|} \hline
 $2.-2\,0.-2.2\,1\,1\,0$ & $2:(-2,2\,1)\,0:-2\,0$ & $2:(2,-2\,-1)\,0:-2\,0$ \\ \hline
 $2.2.-2.2\,0.-2\,-1$ & $3.-2\,0.-2.2\,1\,0$ & $3.-2\,-1\,0.-2.2\,0$ \\ \hline
  $8^*2\,0.-2\,0.-2\,0.2\,0$ & $8^*-2\,-1.2\,0.-2$ & $9^*.-2:-2\,0.-2$ \\ \hline
  $101^*-2\,0::.-2\,0$ & $102^*-2\,0::-2$ &  \\ \hline

\end{tabular}

\normalsize

\bigskip

\noindent while inadequacy of the remaining 8 knots from the
following table remains unknown

\bigskip

\small

\begin{tabular}{|c|c|c|} \hline
$2.-3\,0.-2\,-1.2\,0$ & $8^*-2\,-1\,-1::-2\,0$ &
$8^*2:.-2\,0:.-2\,-1\,0$ \\ \hline

$8^*-2\,-1::-3\,0$ & $8^*2:.-2\,-1\,0:.-2\,0$ &
$8^*-2\,0.2:-2\,-1\,0$ \\ \hline

 $8^*-2.2.-2\,0:2\,0$ & $8^*-2\,0:-2\,0:-2\,0:2\,0$ & \\ \hline

\end{tabular}

\normalsize

\bigskip

For $n=11$ four links

\bigskip

\small

\noindent \begin{tabular}{|c|c|c|c|} \hline

$(2\,1,2)\,-1\,-1\,(2,2)$& $(2,2),-2,-1,(2,-2)$ & $6^*3.-2\,0.-2.2\,0$ & $6^*(2,-2).2.-2$ \\
\hline

\end{tabular}

\bigskip
\normalsize

\noindent are inadequate according to Theorem 6, and the following
8 links are candidates for inadequate links. All their minimal
diagrams are inadequate.

\bigskip

\footnotesize

\noindent \begin{tabular}{|c|c|c|c|} \hline

$(-2\,-1,2)\,1\,1\,(2,-2)$& $(-2\,-1,2,2)\,(2,-2)$ & $(2\,1,2,-2)\,(2,-2)$ & $(2,2)\,-1\,-1\,-1\,(2,-2)$ \\
\hline

$(2,2,-2)\,(-2\,-1,2)$& $(2,2,-2)\,(2\,1,-2)$ &
$6^*2\,1.-2\,0.-2.2\,0$ & $6^*(2,-2),-2$ \\ \hline

\end{tabular}

\bigskip
\normalsize

For $n=12$ inadequacy of 63 links is confirmed according to
Theorem 6, and the remaining 232 links are candidates for
inadequate links.

Tables of adequate non-alternating links with at most $n=12$
crossings in Conway notation can be downloaded in the form of {\it
Mathematica} notebook from the address:

http://www.mi.sanu.ac.yu/vismath/adequate.pdf

\section{Families and classes of links and their adequacy}

\begin{definition}
For a link $L$ given in an unreduced\footnote{The Conway notation
is called unreduced if in symbols of polyhedral links elementary
tangles 1 in single vertices are not omitted.} Conway notation
$C(L)$, let $S$ denote a set of numbers in the Conway symbol,
excluding numbers denoting basic polyhedra and zeros (marking the
position of tangles in the vertices of polyhedra). For $C(L)$ and
an arbitrary (non-empty) subset $\tilde S$ of $S$ the family
$F_{\tilde S}(L)$ of knots or links derived from $L$ is
constructed by substituting each $a \in S_f$, $a\neq 1$, by
$sgn(a) (|a|+k_a)$ for $k_a \in N$ [JaSaz].
\end{definition}

If $k_a$ is an even number ($k_a\in N$), the number of components
is preserved inside a family, i.e., we obtain families of knots or
links with the same number of components.

\begin{definition}
A link given by Conway symbol containing only tangles $1$, $-1$,
$2$, or $-2$ is called a {\it source link}. A link given by Conway
symbol containing only tangles $1$, $-1$, $2$, $-2$, $3$, or $-3$
is called a {\it generating link}.
\end{definition}

\begin{theorem}
All link diagrams which belong to the same family of diagrams have
adequacy of the same sign.
\end{theorem}

{\bf Proof:} If we substitute $a \in S_f$, $a\neq 1$, by $sgn(a)
(|a|+1)$ (Definition 5), a new state circle of the length 2
appears in one of the states $D_{s_+}$ or $D_{s_-}$, so the sign
of adequacy remains unchanged. In the remaining state the number
of state circles remains unchanged and all state circles
associated with the new crossing obtain one new touching point. If
the crossings of the original tangle $a$ after smoothing
correspond to different state circles, the same holds for the
tangle $sgn(a) (|a|+1)$, and the sign of adequacy remains
unchanged. By induction, we conclude that this property holds for
every $k_a \in N$\footnote{See Def. 5.}. Hence, all link diagrams
which belong to the same family of diagrams have the adequacy of
the same sign. $\Box $

{\bf Proposition 1.} The adequacy of a link diagram remains
unchanged if we replace every positive rational tangle by $2$, and
every negative rational tangle by $-2$.

\noindent The proof of this Proposition is straightforward,
because every rational alternating tangle is adequate, so its
collapse into a bigon doe's not change the sign of adequacy.

A pretzel tangle and the pretzel link obtained as its closure,
consisting from $n$ alternating rational tangles $t_i$ is denoted
by $t_1,t_2,\ldots,t_n$ ($n\ge 3$, $t_i\neq 1$, $i=1,...,n$).
Number $n$ will be called the {\it length of the pretzel tangle}.

\begin{theorem}
A non-alternating pretzel link $t_1,t_2,\ldots,t_n$ is
semi-adequate if it contains exactly one rational tangle of one
sign, and all the other rational tangles of the opposite sign.
Otherwise, it is adequate.
\end{theorem}

A pretzel tangle is called adequate or semi-adequate if its
corresponding pretzel link is adequate or semi-adequate,
respectively.

Let's denote source link of the form $2,...,-2,...$, where 2
occurs $k$ times, and $-2$ occurs $l$ times with the short symbol
$(2)^k,(-2)^l$. For different values of $k$ and $l$ we obtain six
classes of source links, where all members of the same class have
the adequacy of the same sign

\bigskip

\begin{tabular}{|c|c|} \hline
$k\ge 3$, $l=0$   & $+$alternating \\ \hline

$k=0$, $l\ge 3$   & $-$alternating \\ \hline

$k\ge 2$, $l\ge 2$   &  adequate \\ \hline

$k=1$, $l\ge 2$   & $+$adequate \\ \hline

$k\ge 2$, $l=1$  & $-$adequate \\ \hline

$k=l=1$ & inadequate \\ \hline
\end{tabular}

\bigskip
This property directly follows from Theorem 1 and Theorem 8. As
the minimal representatives of these six classes we can use source
links $(2,2,2)$, $(-2,-2,-2)$, $(2,2,-2,-2)$, $(-2,-2,2)$,
$(2,2,-2)$, and $(2,-2)$, respectively. Combining this with
Proposition 1 we conclude that these six source links can be used
as the representatives of the corresponding pretzel links
(Montesinos links) with the rational tangles of the corresponding
signs. For example, source link $2,2,-2,-2$ can be used as the
representative of all non-alternating adequate pretzel links of
the form $t_1,...,t_k,-t'_1,...,-t'_l$, ($k\ge 2$, $l\ge 2$),
where $t_i$ ($i=1,2,...k$) and $t'_j$ ($j=1,2,...,l$) are positive
rational tangles different from 1.

\section{Some particular classes of algebraic links and their adequacy}

\begin{definition}
An alternating pretzel tangle $P_n=t_1,t_2,\ldots,t_n$ is called
$+$alternating if all its rational tangles $t_i$ are positive, and
$-$alternating if they are all negative.
\end{definition}

Tangle $t_1,-t_2$ is inadequate, where $t_1$, $t_2$ are positive
rational tangles.

\begin{theorem}
A link $P_m\,Q_n=(p_1,p_2,\ldots,p_m)\,(q_1,q_2,\ldots,q_n)$
($m,n\ge 2$) obtained as the product of pretzel tangles $P_m$ and
$Q_n$ is adequate if
\begin{itemize}
\item both $P_m$ and $Q_n$ are adequate; or
\item one of them is $+$alternating, and the other $+$adequate; or
\item one of them is $-$alternating, and the other $-$adequate.
\end{itemize}
It is semi-adequate if
\begin{itemize}
\item one of them is adequate, and the other semi-adequate; or
\item one of them is $+$adequate, and the other $-$adequate; or
\item if one of them inadequate, and the other an
alternating pretzel tangle.
\end{itemize}
It is candidate for inadequate if
\begin{itemize}
\item both $P_m$ and $Q_n$ are $+$adequate or $-$adequate;
\item if one of them is inadequate, and the other is not an
alternating pretzel tangle.
\end{itemize}
\end{theorem}

From the preceding theorem we obtain the following multiplication
table, where the * denotes the product of pretzel
tangles\footnote{The product $P_1\,P_2$ of inadequate tangles
$P_1$ and $P_2$ is omitted, since it represents a non-minimal
diagram of an alternating link.}:

\bigskip

\begin{tabular}{|c|c|c|c|c|c|c|} \hline
  {\bf *} & {\bf $+$alt}  & {\bf $-$alt}  & {\bf adq} & {\bf $+$adq} & {\bf $-$adq} & {\bf inadq} \\ \hline

 {\bf $+$alt}  & $+$alt & adq & adq & adq & $+$adq & $+$adq \\ \hline

 {\bf $-$alt}  & adq & $-$alt & adq & $-$adq & adq & $-$adq \\ \hline

 {\bf adq}  & adq & adq & adq & $-$adq & $+$adq  & inadeq \\ \hline

 {\bf $+$adq}  & adq & $+$adq & $+$adq & inadeq & $+$adq & inadeq \\ \hline

 {\bf $-$adq}  & $-$adq & adq & $-$adq & $-$adq & inadeq & inadeq \\ \hline

 {\bf inadq}  & $-$adq & $+$adq & inadeq  & inadeq & inadeq &  \\ \hline

\end{tabular}

\bigskip

For links of the form
$P_m\,Q_n=(p_1,p_2,\ldots,p_n)\,(q_1,q_2,\ldots,q_n)$ we obtained
general rules for adequacy, expressed as the multiplication table.
Unfortunately, for links of the form $P_1\,P_2\,...\,P_k$, with
$k\ge 2$ we are not able to present general adequacy
multiplication tables.

As the minimal representatives of pretzel tangles with the
properties {\bf $+$alt}, {\bf $-$alt}, {\bf adq}, {\bf $+$adq},
{\bf $-$adq}, and {\bf inadeq} we can use the following tangles:

\bigskip

\begin{tabular}{|c|c|c|} \hline
 1  & {\bf $+$alt} & $2,2,2$\\ \hline

 2  & {\bf $-$alt} & $-2,-2,-2$ \\ \hline

 3  & {\bf adq} & $2,2,-2,-2$ \\ \hline

 4  & {\bf $+$adq} & $-2,-2,2$ \\ \hline

 5  & {\bf $-$adq} & $2,2,-2$ \\ \hline

 6  & {\bf inadeq} & $2,-2$ \\ \hline
\end{tabular}

\bigskip

If we denote the properties {\bf $+$alt}, {\bf $-$alt}, {\bf adq},
{\bf $+$adq}, {\bf $-$adq}, and {\bf inadeq} by 1-6, for $k=3$, we
have the following statement:

\begin{theorem}
The links $P_1\,P_2\,P_3$ are adequate for the following
properties of pretzel tangles $P_1$, $P_2$, $P_3$:

\small

\bigskip
\begin{tabular}{|c|c|c|c|c|c|c|c|} \hline
 $1,1,1$  & $1,1,2$ & $1,1,3$ & $1,1,4$ & $1,2,1$ & $1,2,2$ & $1,2,3$ &  $1,3,1$  \\ \hline

 $1,3,2$  & $1,3,3$ & $1,4,1$ & $1,4,2$ & $1,4,3$ & $1,5,2$ & $1,6,2$ &  $2,1,2$  \\ \hline

 $2,1,3$  & $2,2,2$ & $2,2,3$ & $2,2,5$ & $2,3,2$ & $2,3,3$ & $2,5,2$ &  $2,5,3$  \\ \hline

 $3,1,3$  & $3,1,4$ & $3,2,3$ & $3,2,5$ & $3,3,3$ & $4,1,4$ & $5,2,5$ &  $$  \\ \hline

\end{tabular}

\normalsize

\bigskip

\noindent semi-adequate for:

\small

\bigskip
\begin{tabular}{|c|c|c|c|c|c|c|c|} \hline

 $1,1,5$  & $1,1,6$ & $1,2,4$ & $1,2,5$ & $1,3,4$ & $1,3,5$ & $1,4,4$ &  $1,4,5$  \\ \hline

 $1,5,1$  & $1,5,3$ & $1,5,4$ & $1,6,1$ & $1,6,3$ & $1,6,4$ & $2,1,4$ &  $2,1,5$  \\ \hline

 $2,2,4$  & $2,2,6$ & $2,3,4$ & $2,3,5$ & $2,4,2$ & $2,4,3$ & $2,4,5$ &  $2,5,4$  \\ \hline

 $2,5,5$  & $2,6,2$ & $2,6,3$ & $2,6,5$ & $3,1,5$ & $3,1,6$ & $3,2,4$ &  $3,2,6$  \\ \hline

 $3,3,4$  & $3,3,5$ & $3,4,3$ & $3,4,5$ & $3,5,3$ & $3,5,4$ & $4,1,5$ &  $4,1,6$  \\ \hline

 $4,2,4$  & $4,2,5$ & $4,2,6$ & $4,3,4$ & $4,5,4$ & $5,1,5$ & $5,1,6$ &  $5,2,6$  \\ \hline

 $5,3,5$  & $5,4,5$ & $6,1,6$ & $6,2,6$ & $$ & $$ & $$ &  $$  \\ \hline

\end{tabular}

\normalsize

\bigskip

\noindent and candidates for inadequate for:

\small

\bigskip
\begin{tabular}{|c|c|c|c|c|c|c|c|} \hline

  $1,2,6$  & $1,3,6$ & $1,4,6$ & $1,5,5$ & $1,5,6$ & $1,6,5$ & $1,6,6$ &  $2,1,6$  \\ \hline

  $2,3,6$  & $2,4,4$ & $2,4,6$ & $2,5,6$ & $2,6,4$ & $2,6,6$ & $3,3,6$ &  $3,4,4$  \\ \hline

  $3,4,6$  & $3,5,5$ & $3,5,6$ & $3,6,3$ & $3,6,4$ & $3,6,5$ & $3,6,6$ &  $4,3,5$  \\ \hline

  $4,3,6$  & $4,4,4$ & $4,4,5$ & $4,4,6$ & $4,5,5$ & $4,5,6$ & $4,6,4$ &  $4,6,5$  \\ \hline

  $4,6,6$  & $5,3,6$ & $5,4,6$ & $5,5,5$ & $5,5,6$ & $5,6,5$ & $5,6,6$ &  $6,3,6$  \\ \hline

  $6,4,6$  & $6,5,6$ & $6,6,6$ & $$ & $$ & $$ & $$ &  $$  \\ \hline

\end{tabular}

\normalsize

\bigskip

\end{theorem}

The results hold for all sequences $a,b,c$ ($a,b,c\in
\{1,2,...,6\}$) and their reverses. Analogous tables are obtained
by computer calculations for all $k\le 6$.

For a given non-alternating pretzel tangle $P$ the tangle $P'$
obtained from it by replacing every rational positive or negative
tangle $t_i$ with the tangle $sign(t_i)\times 2$ will be called
basic pretzel tangle.

\begin{theorem}
The links $P_1\,P_2\,...\,P_k$ and $P_1'\,P_2'\,...\,P_k'$ have
the same adequacy.
\end{theorem}

Next, we will consider links of the form $P_1,P_2,...,P_k$ ($k\ge
3$, where $P_i$ ($i=1,...,k$) are pretzel tangles. Since
permutation of pretzel tangles preserves the sign of adequacy, the
result holds for every sequence $a,b,c$ ($a,b,c\in \{1,2,...,6\}$)
and all of its permutations. For $k=3$ we obtained the following
result:

\begin{theorem}
The links $P_1,P_2,P_3$ are adequate for the following properties
of pretzel tangles $P_1$, $P_2$, $P_3$:

\small

\bigskip
\begin{tabular}{|c|c|c|c|c|c|c|c|} \hline
 $1,1,1$  & $1,1,2$ & $1,1,3$ & $1,1,4$ & $1,2,2$ & $1,2,3$ & $1,3,3$ &  $1,3,4$  \\ \hline

 $1,4,4$  & $2,2,2$ & $2,2,3$ & $2,2,5$ & $2,3,3$ & $2,3,5$ & $2,5,5$ &  $3,3,3$  \\ \hline

 $3,3,4$  & $3,3,5$ & $3,3,6$ & $3,4,4$ & $3,4,5$ & $3,4,6$ & $3,5,5$ &  $3,5,6$  \\ \hline

 $3,6,6$  & $4,4,4$ & $4,4,5$ & $4,4,6$ & $4,5,5$ & $4,5,6$ & $4,6,6$ &  $5,5,5$  \\ \hline

  $5,5,6$  & $5,6,6$ & $6,6,6$ & $$ & $$ & $$ & $$ &  $$  \\ \hline

\end{tabular}

\normalsize

\bigskip

\noindent semi-adequate for:

\small

\bigskip
\begin{tabular}{|c|c|c|c|c|c|c|c|} \hline

   $1,2,4$  & $2,2,4$ & $2,2,6$ & $2,3,4$ & $2,3,6$ & $2,4,4$ & $2,4,5$ &  $2,4,6$  \\ \hline

   $2,5,6$  & $2,6,6$ & $1,1,5$ & $1,1,6$ & $1,2,5$ & $1,3,5$ & $1,3,6$ &  $1,4,5$  \\ \hline

   $1,4,6$  & $1,5,5$ & $1,5,6$ & $1,6,6$ & $$ & $$ & $$ &  $$  \\ \hline

\end{tabular}

\normalsize

\bigskip

\noindent and candidates for inadequate for:

\small

\bigskip
\begin{tabular}{|c|} \hline

  $1,2,6$  \\ \hline

\end{tabular}

\normalsize

\bigskip

\end{theorem}

\noindent Analogous results are obtained by computer calculations
for all $k\le 6$.

Next we consider links of the form
$(P_1,P_2,...,P_m)\,(Q_1,Q_2,...,Q_n)$ ($m,n\ge 2$), where $P_i$
and $Q_j$ ($i=1,...,m$, $j=1,...,n$) are pretzel tangles.

In the case $m=n=2$, where in the sequence $(a,b)\,(c,d)$
($a,b,c,d\in \{1,2,...,6\}$, $a$ and $b$, and $c$ and $d$ can
commute and all sequences $(a,b)\,(c,d)$ can be reversed, we
obtain the following result:

\begin{theorem}
The links $(P_1,P_2)\,(Q_1,Q_2)$ are adequate for the following
properties of pretzel tangles $P_1$, $P_2$, $Q_1$, $Q_2$:

\tiny

\bigskip
\noindent \begin{tabular}{|c|c|c|c|c|c|c|} \hline
 $(1,1)\,(1,1)$  & $(1,1)\,(1,2)$ & $(1,1)\,(1,3)$ & $(1,1)\,(1,4)$ & $(1,1)\,(2,2)$ & $(1,1)\,(2,3)$ & $(1,1)\,(3,3)$   \\ \hline
 $(1,1)\,(3,4)$ &

 $(1,1)\,(4,4)$  & $(1,2)\,(1,2)$ & $(1,2)\,(1,3)$ & $(1,2)\,(1,4)$ & $(1,2)\,(2,2)$ & $(1,2)\,(2,3)$ \\ \hline

  $(1,2)\,(2,5)$ &  $(1,2)\,(3,3)$ &

 $(1,2)\,(3,4)$  & $(1,2)\,(3,5)$ & $(1,2)\,(3,6)$ & $(1,2)\,(4,4)$ & $(1,2)\,(4,5)$ \\ \hline

 $(1,2)\,(4,6)$ & $(1,2)\,(5,5)$ &  $(1,2)\,(5,6)$ &

 $(1,2)\,(6,6)$  & $(1,3)\,(1,3)$ & $(1,3)\,(1,4)$ & $(1,3)\,(2,2)$ \\ \hline

  $(1,3)\,(2,3)$ & $(1,3)\,(3,3)$ & $(1,3)\,(3,4)$ &
  $(1,3)\,(4,4)$ &   $(1,4)\,(1,4)$  & $(2,2)\,(2,2)$ & $(2,2)\,(2,3)$ \\ \hline

  $(2,2)\,(2,5)$ & $(2,2)\,(3,3)$ & $(2,2)\,(3,5)$ & $(2,2)\,(5,5)$ &
  $(2,3)\,(2,3)$ &   $(2,3)\,(2,5)$ &

   $(2,3)\,(3,3)$ \\ \hline  $(2,3)\,(3,5)$ & $(2,3)\,(5,5)$ & $(2,5)\,(2,5)$ & $(3,3)\,(3,3)$ & $$ &  $$ & \\ \hline

\end{tabular}

\normalsize

\bigskip

\noindent semi-adequate for:

\tiny

\bigskip
\noindent\begin{tabular}{|c|c|c|c|c|c|c|} \hline

   $(1,1)\,(1,5)$  & $(1,1)\,(1,6)$ & $(1,1)\,(2,5)$ & $(1,1)\,(3,5)$ & $(1,1)\,(3,6)$ & $(1,1)\,(4,5)$ & $(1,1)\,(4,6)$   \\ \hline

   $(1,1)\,(5,5)$  & $(1,1)\,(5,6)$ & $(1,1)\,(6,6)$ & $(1,2)\,(1,5)$ & $(1,2)\,(1,6)$ & $(1,3)\,(1,5)$ & $(1,3)\,(1,6)$ \\ \hline

   $(1,3)\,(2,5)$  & $(1,3)\,(3,5)$ & $(1,3)\,(3,6)$ & $(1,3)\,(4,5)$ & $(1,3)\,(4,6)$ & $(1,3)\,(5,5)$ & $(1,3)\,(5,6)$  \\ \hline

   $(1,3)\,(6,6)$  & $(1,4)\,(1,5)$ & $(1,4)\,(1,6)$ & $(1,4)\,(2,2)$ & $(1,4)\,(2,3)$ & $(1,4)\,(2,5)$ & $(1,4)\,(3,3)$  \\ \hline

   $(1,4)\,(3,4)$  & $(1,4)\,(3,5)$ & $(1,4)\,(3,6)$ & $(1,4)\,(4,4)$ & $(1,4)\,(4,5)$ & $(1,4)\,(4,6)$ & $(1,4)\,(5,5)$  \\ \hline

   $(1,4)\,(5,6)$  & $(1,4)\,(6,6)$ & $(2,4)\,(2,5)$ & $(2,4)\,(3,3)$ & $(2,4)\,(3,5)$ & $(2,4)\,(5,5)$ & $(3,3)\,(3,5)$  \\ \hline

   $(3,3)\,(5,5)$  & $(3,4)\,(3,5)$ & $(3,4)\,(5,5)$ & $(4,4)\,(5,5)$ & $(1,1)\,(2,4)$ & $(1,2)\,(2,4)$ & $(1,2)\,(2,6)$  \\ \hline

   $(1,3)\,(2,4)$  & $(1,5)\,(2,2)$ & $(1,5)\,(2,3)$ & $(1,5)\,(2,4)$ & $(1,5)\,(3,3)$ & $(1,5)\,(3,4)$ & $(1,5)\,(4,4)$  \\ \hline

   $(2,2)\,(2,4)$  & $(2,2)\,(2,6)$ & $(2,2)\,(3,4)$ & $(2,2)\,(3,6)$ & $(2,2)\,(4,4)$ & $(2,2)\,(4,5)$ & $(2,2)\,(4,6)$  \\ \hline

   $(2,2)\,(5,6)$  & $(2,2)\,(6,6)$ & $(2,3)\,(2,4)$ & $(2,3)\,(2,6)$ & $(2,3)\,(3,4)$ & $(2,3)\,(3,6)$ & $(2,3)\,(4,4)$  \\ \hline

   $(2,3)\,(4,5)$  & $(2,3)\,(4,6)$ & $(2,3)\,(5,6)$ & $(2,3)\,(6,6)$ & $(2,5)\,(2,6)$ & $(2,5)\,(3,3)$ & $(2,5)\,(3,4)$  \\ \hline

   $(2,5)\,(3,5)$  & $(2,5)\,(3,6)$ & $(2,5)\,(4,4)$ & $(2,5)\,(4,5)$ & $(2,5)\,(4,6)$ & $(2,5)\,(5,5)$ & $(2,5)\,(5,6)$  \\ \hline

   $(2,5)\,(6,6)$  & $(3,3)\,(3,4)$ & $(3,3)\,(4,4)$ & $(3,5)\,(4,4)$ & $$ & $$ & $$  \\ \hline

\end{tabular}

\normalsize

\bigskip

\noindent and candidates for inadequate for:

\tiny

\bigskip
\noindent\begin{tabular}{|c|c|c|c|c|c|c|} \hline

   $(1,1)\,(2,6)$  & $(1,3)\,(2,6)$ & $(1,4)\,(2,4)$ & $(1,4)\,(2,6)$ & $(1,5)\,(1,5)$ & $(1,5)\,(1,6)$ & $(1,5)\,(2,5)$  \\ \hline

   $(1,5)\,(2,6)$  & $(1,5)\,(3,5)$ & $(1,5)\,(3,6)$ & $(1,5)\,(4,5)$ & $(1,5)\,(4,6)$ & $(1,5)\,(5,5)$ & $(1,5)\,(5,6)$  \\ \hline

   $(1,5)\,(6,6)$  & $(1,6)\,(1,6)$ & $(1,6)\,(2,2)$ & $(1,6)\,(2,3)$ & $(1,6)\,(2,4)$ & $(1,6)\,(2,5)$ & $(1,6)\,(2,6)$  \\ \hline

   $(1,6)\,(3,3)$  & $(1,6)\,(3,4)$ & $(1,6)\,(3,5)$ & $(1,6)\,(3,6)$ & $(1,6)\,(4,4)$ & $(1,6)\,(4,5)$ & $(1,6)\,(4,6)$  \\ \hline

   $(1,6)\,(5,5)$  & $(1,6)\,(5,6)$ & $(1,6)\,(6,6)$ & $(2,4)\,(2,4)$ & $(2,4)\,(2,6)$ & $(2,4)\,(3,4)$ & $(2,4)\,(3,6)$  \\ \hline

   $(2,4)\,(4,4)$  & $(2,4)\,(4,5)$ & $(2,4)\,(4,6)$ & $(2,4)\,(5,6)$ & $(2,4)\,(6,6)$ & $(2,6)\,(2,6)$ & $(2,6)\,(3,3)$  \\ \hline

   $(2,6)\,(3,4)$  & $(2,6)\,(3,5)$ & $(2,6)\,(3,6)$ & $(2,6)\,(4,4)$ & $(2,6)\,(4,5)$ & $(2,6)\,(4,6)$ & $(2,6)\,(5,5)$  \\ \hline

   $(2,6)\,(5,6)$  &  $(2,6)\,(6,6)$ & $(3,3)\,(3,6)$ & $(3,3)\,(4,5)$ & $(3,3)\,(4,6)$ & $(3,3)\,(5,6)$ & $(3,3)\,(6,6)$  \\ \hline

   $(3,4)\,(3,4)$ & $(3,4)\,(3,6)$ & $(3,4)\,(4,4)$ & $(3,4)\,(4,5)$ & $(3,4)\,(4,6)$ & $(3,4)\,(5,6)$ & $(3,4)\,(6,6)$  \\ \hline

   $(3,5)\,(3,5)$  & $(3,5)\,(3,6)$ & $(3,5)\,(4,5)$ & $(3,5)\,(4,6)$ & $(3,5)\,(5,5)$ & $(3,5)\,(5,6)$ & $(3,5)\,(6,6)$  \\ \hline

   $(3,6)\,(3,6)$  & $(3,6)\,(4,4)$ & $(3,6)\,(4,5)$ & $(3,6)\,(4,6)$ & $(3,6)\,(5,5)$ & $(3,6)\,(5,6)$ & $(3,6)\,(6,6)$  \\ \hline

   $(4,4)\,(4,4)$  & $(4,4)\,(4,5)$ & $(4,4)\,(4,6)$ & $(4,4)\,(5,6)$ & $(4,4)\,(6,6)$ & $(4,5)\,(4,5)$ & $(4,5)\,(4,6)$  \\ \hline

   $(4,5)\,(5,5)$  & $(4,5)\,(5,6)$ & $(4,5)\,(6,6)$ & $(4,6)\,(4,6)$ & $(4,6)\,(5,5)$ & $(4,6)\,(5,6)$ & $(4,6)\,(6,6)$  \\ \hline

   $(5,5)\,(5,5)$  & $(5,5)\,(5,6)$ & $(5,5)\,(6,6)$ & $(5,6)\,(5,6)$ & $(5,6)\,(6,6)$ & $(6,6)\,(6,6)$ & $$  \\ \hline

\end{tabular}

\normalsize

\bigskip

\end{theorem}

\noindent Analogous results are obtained by computer calculations
for $m,n\le 4$.

Furthermore we consider links of the form $P_1,t_1,t_2,...,t_n$,
where $P_1$ is a pretzel tangle, and $t_i$ ($i=1,2,...n$, $n\ge
2$) are rational tangles. If $P=t_1,t_2,...,t_n$, we have the
following statement:

\begin{itemize}

\item links of the given form are adequate if $\{P_1,P\}\in
\{\{1,1\},\{1,2\},\{1,3\},\{1,4\},$
$\{2,2\},\{2,3\},\{2,5\},\{3,3\}\}$;

\item semi-adequate if $\{P_1,P\}\in
\{\{1,5\},\{1,6\},\{2,4\},\{2,6\},\{3,4\},$ $\{3,5\},\{4,5\}\}$;

\item and candidates for inadequate if $\{P_1,P\}\in
\{\{3,6\},\{4,4\},\{4,6\},\{5,5\},\{5,6\},$ $\{6,6\}\}$.

\end{itemize}

Next we consider links of the form $P_1,...,P_m,t_1,t_2,...,t_n$,
where $P_1$,...,$P_m$ are pretzel tangles, and $t_i$
($i=1,2,...n$, $n\ge 2$) are rational tangles. If
$P=t_1,t_2,...,t_n$, for $m=2$ we have the following statement:

Links of the given form are adequate if $(\{P_1,P_2\},P)$ is:

\small

\bigskip
\begin{tabular}{|c|c|c|c|c|c|} \hline

   $(\{1,1\},1)$  & $(\{1,1\},2)$ & $(\{1,1\},3)$ & $(\{1,1\},4)$ & $(\{1,2\},1)$ & $(\{1,2\},2)$   \\ \hline

   $(\{1,2\},3)$  & $(\{1,2\},4)$ & $(\{1,2\},5)$ & $(\{1,2\},6)$ & $(\{1,3\},1)$ & $(\{1,3\},2)$   \\ \hline

   $(\{1,3\},3)$  & $(\{1,3\},4)$ & $(\{1,4\},1)$ & $(\{2,2\},1)$ & $(\{2,2\},2)$ & $(\{2,2\},3)$  \\ \hline

   $(\{2,2\},5)$  & $(\{2,3\},1)$ & $(\{2,3\},2)$ & $(\{2,3\},3)$ & $(\{2,3\},5)$ & $(\{2,5\},2)$  \\ \hline

   $(\{3,3\},1)$  & $(\{3,3\},2)$ & $(\{3,3\},3)$ & $(\{3,4\},1)$ & $(\{3,5\},2)$ & $(\{4,4\},1)$   \\ \hline

   $(\{5,5\},2)$  & $$ & $$ & $$ & $$ & $$   \\ \hline

\end{tabular}

\normalsize

\bigskip

semi-adequate if $(\{P_1,P_2\},P)$ is:

\small

\bigskip
\begin{tabular}{|c|c|c|c|c|c|} \hline

   $(\{1,1\},5)$  & $(\{1,1\},6)$ & $(\{1,3\},5)$ & $(\{1,3\},6)$ & $(\{1,4\},2)$ & $(\{1,4\},3)$   \\ \hline

   $(\{1,4\},4)$  & $(\{1,4\},5)$ & $(\{1,4\},6)$ & $(\{1,5\},1)$ & $(\{1,5\},2)$ & $(\{1,5\},3)$   \\ \hline

   $(\{1,5\},4)$  & $(\{1,6\},1)$ & $(\{2,2\},4)$ & $(\{2,2\},6)$ & $(\{2,3\},4)$ & $(\{2,3\},6)$   \\ \hline

   $(\{2,4\},1)$  & $(\{2,4\},2)$ & $(\{2,4\},3)$ & $(\{2,4\},5)$ & $(\{2,5\},1)$ & $(\{2,5\},3)$   \\ \hline

   $(\{2,5\},4)$  & $(\{2,5\},5)$ & $(\{2,5\},6)$ & $(\{2,6\},2)$ & $(\{3,3\},4)$ & $(\{3,3\},5)$   \\ \hline

   $(\{3,4\},2)$  & $(\{3,4\},3)$ & $(\{3,4\},5)$ & $(\{3,5\},1)$ & $(\{3,5\},3)$ & $(\{3,5\},4)$   \\ \hline

   $(\{3,6\},1)$  & $(\{3,6\},2)$ & $(\{4,4\},2)$ & $(\{4,4\},3)$ & $(\{4,4\},5)$ & $(\{4,5\},1)$   \\ \hline

   $(\{4,5\},2)$  & $(\{4,6\},1)$ & $(\{4,6\},2)$ & $(\{5,5\},1)$ & $(\{5,5\},3)$ & $(\{5,5\},4)$   \\ \hline

   $(\{5,6\},1)$  & $(\{5,6\},2)$ & $(\{6,6\},1)$ & $(\{6,6\},2)$ & $$ & $$   \\ \hline

\end{tabular}

\normalsize

\bigskip

and candidates for inadequate if $(\{P_1,P_2\},P)$ is:

\small

\bigskip
\begin{tabular}{|c|c|c|c|c|c|} \hline

   $(\{1,5\},5)$  & $(\{1,5\},6)$ & $(\{1,6\},2)$ & $(\{1,6\},3)$ & $(\{1,6\},4)$ & $(\{1,6\},5)$   \\ \hline

   $(\{1,6\},6)$  & $(\{2,4\},4)$ & $(\{2,4\},6)$ & $(\{2,6\},1)$ & $(\{2,6\},3)$ & $(\{2,6\},4)$   \\ \hline

   $(\{2,6\},5)$  & $(\{2,6\},6)$ & $(\{3,3\},6)$ & $(\{3,4\},4)$ & $(\{3,4\},6)$ & $(\{3,5\},5)$   \\ \hline

   $(\{3,5\},6)$  & $(\{3,6\},3)$ & $(\{3,6\},4)$ & $(\{3,6\},5)$ & $(\{3,6\},6)$ & $(\{4,4\},4)$   \\ \hline

   $(\{4,4\},6)$  & $(\{4,5\},3)$ & $(\{4,5\},4)$ & $(\{4,5\},5)$ & $(\{4,5\},6)$ & $(\{4,6\},3)$   \\ \hline

   $(\{4,6\},4)$  & $(\{4,6\},5)$ & $(\{4,6\},6)$ & $(\{5,5\},5)$ & $(\{5,5\},6)$ & $(\{5,6\},3)$   \\ \hline

   $(\{5,6\},4)$  & $(\{5,6\},5)$ & $(\{5,6\},6)$ & $(\{6,6\},3)$ & $(\{6,6\},4)$ & $(\{6,6\},5)$   \\ \hline

\end{tabular}

\normalsize

\bigskip

In the same way, by experimental computer calculations, it is
possible to obtain the results for some other types of links. For
example, a link of the form $P_1\,p\,P_2$, where $P_1$, $P_2$ are
pretzel tangles, and $p$ is a positive chain of bigons is adequate
if $\{P_1,P_2\}\in
\{\{1,1\},\{1,3\},\{1,4\},\{3,3\},\{3,4\},\{4,4\}\}$, candidate
for inadequate if $\{P_1,P_2\}\in \{\{2,5\},\{2,6\}\}$, and
semi-adequate otherwise.

\section{Adequacy of polyhedral links}

\begin{theorem}
Polyhedral link with one pretzel tangle $P_1$ and positive
rational tangles in other vertices is adequate if $P_1\in
\{1,2,3,4\}$, and semi-adequate if $P_1\in \{5,6\}$.
\end{theorem}

\begin{theorem}
In every adequate polyhedral link with two pretzel tangles $P_1$,
$P_2$ and positive rational tangles in other vertices, $P_1\notin
\{5,6\}$ and $P_2\notin \{5,6\}$.
\end{theorem}

Condition from the Theorem 15 is necessary, but not sufficient.
Hence, we will consider different cases, depending on different
polyhedral source links. For example, the following results hold for
non-alternating links derived from the basic polyhedron $6^*$ with
two pretzel tangles $P_1$ and $P_2$ and positive rational tangles in
remaining vertices:

\begin{itemize}
\item a link of the form $6^*P_1.P_2.t_1.t_2.t_3.t_4$ is adequate
if $\{P_1,P_2\}\in \{\{1,1\},\{1,2\}\{1,3\},$
$\{1,4\},\{3,3\},\{3,4\},\{4,4\}\}$, a candidate for inadequate if
$\{P_1,P_2\}\in \{\{2,2\},\{5,5\},$ $\{5,6\},\{6,6\}\}$, and
semi-adequate otherwise;
\item a link of the form $6^*P_1.P_2\,0.t_1.t_2.t_3.t_4$ is adequate
if $\{P_1,P_2\}\in \{\{1,1\},\{1,2\},$ $\{1,3\},\{1,4\}\}$, a
candidate for inadequate if $\{P_1,P_2\}\in \{\{2,2\}\}$, and
semi-adequate otherwise;
\item a link of the form $6^*P_1.t_1.t_2.P_2\,0.t_3.t_4$ is adequate
if $P_1\notin \{5,6\}$ and $P_2\notin \{5,6\}$, and semi-adequate
otherwise;
\item a link of the form $6^*P_1.t_1.t_2.P_2.t_3.t_4$ is adequate
if $P_1\notin \{5,6\}$ and $P_2\notin \{5,6\}$, and a candidate
for inadequate if $\{P_1,P_2\}= \{5,6\}$.
\end{itemize}

A link of the form $6^*P_1.P_2.P_3.t_1.t_2.t_3$ is adequate for
the following triples $(P_1,P_2,P_3)$:

\small

\bigskip
\begin{tabular}{|c|c|c|c|c|c|c|c|} \hline
 $1,1,1$  & $1,1,2$ & $1,1,3$ & $1,1,4$ & $1,2,1$ & $1,2,3$ & $1,2,4$ &  $1,3,1$  \\ \hline

 $1,3,3$  & $1,3,4$ & $1,4,1$ & $1,4,3$ & $1,4,4$ & $2,1,1$ & $2,1,3$ &  $2,1,4$  \\ \hline

 $3,1,1$  & $3,1,2$ & $3,1,3$ & $3,1,4$ & $3,2,1$ & $3,2,3$ & $3,2,4$ &  $3,2,5$  \\ \hline

 $3,2,6$  & $3,3,1$ & $3,3,3$ & $3,3,4$ & $3,4,1$ & $4,1,1$ & $4,1,2$ &  $4,1,3$  \\ \hline

 $4,1,4$  & $4,2,1$ & 4,2,3$$ & $4,2,4$ & $4,2,5$ & $4,2,6$ & $4,3,1$ &  $4,3,3$  \\ \hline

 $4,3,4$  & $4,4,1$ & $5,2,3$ & $5,2,4$ & $5,2,5$ & $5,2,6$ & $6,2,3$ &  $6,2,4$  \\ \hline

 $6,2,5$  & $6,2,6$ & $$ & $$ & $$ & $$ & $$ &  $$  \\ \hline

\end{tabular}

\normalsize

\bigskip

\noindent a candidate for inadequate for:

\small

\bigskip
\begin{tabular}{|c|c|c|c|c|c|c|c|} \hline
 $1,5,5$  & $1,5,6$ & $1,6,5$ & $1,6,6$ & $2,3,5$ & $2,3,6$ & $2,4,2$ &  $2,4,3$  \\ \hline

 $2,4,4$  & $2,4,5$ & $2,4,6$ & $2,5,5$ & $2,5,6$ & 2,6,2$$ & $2,6,3$ &  $2,6,4$  \\ \hline

 $2,6,5$  & $2,6,6$ & $3,4,2$ & $3,5,5$ & $3,5,6$ & $3,6,2$ & $3,6,3$ &  $3,6,4$  \\ \hline

 $3,6,5$  & $3,6,6$ & $4,4,2$ & $4,5,5$ & $4,5,6$ & $4,6,2$ & $4,6,3$ &  $4,6,4$  \\ \hline

 $4,6,5$  & $4,6,6$ & $5,3,2$ & $5,4,2$ & $5,5,1$ & $5,5,2$ & $5,5,3$ &  $5,5,4$  \\ \hline

 $5,5,5$  & $5,5,6$ & 5,6,1$$ & $5,6,2$ & $5,6,3$ & $5,6,4$ & $5,6,5$ &  $5,6,6$  \\ \hline

 $6,3,2$  & $6,4,2$ & $6,5,1$ & $6,5,2$ & $6,5,3$ & $6,5,4$ & $6,5,5$ &  $6,5,6$  \\ \hline

 $6,6,1$  & $6,6,2$ & $6,6,3$ & $6,6,4$ & $6,6,5$ & $6,6,6$ & $$ &  $$  \\ \hline

\end{tabular}

\normalsize

\bigskip

\noindent and semi-adequate otherwise.

Except polyhedral links with pretzel tangles, we will consider
polyhedral links containing only rational tangles.

\begin{theorem}
Non-alternating link derived from the basic polyhedron $6^*$ is a
candidate for inadequate if it is obtained from one of the
following source links by replacing $2$-tangles by positive
rational tangles $t_i$ ($i\in \{1,...,6\}$, $t_i\neq 1$)

\small

\bigskip
\noindent\begin{tabular}{|c|c|c|} \hline

$6^*2.-2\,0.-2.2\,0$  &  $6^*2.2.-2.2.-2\,0$  &
$6^*-2.2.-2\,0.2.2$   \\ \hline

$6^*2.2.-2.2\,0.-2$  &  $6^*-2.2\,0.-2.2\,0.2$  &
$6^*2.-2.2.2.2\,0.-2\,0$   \\ \hline

$6^*2.-2.2\,0.-2.-2.-2\,0$  &  $6^*2.-2.-2.-2.2.-2\,0$  &  $6^*2.-2.2.2.2.-2\,0$   \\
\hline

$6^*2.-2.-2.-2\,0.2.-2\,0$  &  $6^*2.-2\,0.-2.-2\,0.-2.2\,0$  &
$6^*2.-2\,0.-2.2\,0.-2.2\,0$
\\ \hline

\end{tabular}

\normalsize

\bigskip

\noindent and semi-adequate otherwise\footnote{Knot
$6^*2.-2.2\,0.-2.-2.-2\,0$ is recognized as potential inadequate,
i.e., as a knot without minimal $+$ or $-$adequate diagram by
M.~Thistlethwaite in 1988 [Thi], but to this knot cannot be appied
Theorem 6. From 12 links from this table, the five of them:
$6^*2.-2\,0.-2.2\,0$, $6^*2.2.-2.2.-2\,0$, $6^*2.2.-2.2\,0.-2$,
$6^*2.-2.2.2.2\,0.-2\,0$, and $6^*2.-2.-2.-2.2.-2\,0$ are
inadequate, according to Theorem 6.}.

\end{theorem}

\noindent In the same way, similar results is possible to obtain
for other basic polyhedra.

\section{Adequacy of mixed states and adequacy number}

The definition of adequacy can be extended to an arbitrary state a
link diagram $D$. Together with special states $s_+$ and $s_-$, we
will consider mixed states, where markers have different signs.

According to Definition 1, a state $s$ of the diagram $D$ is
called {\it adequate state} if, at each crossing, the two segments
of $D_s$ which replace the crossing are in different state
circles.

\begin{theorem}
Every  link diagram has at least two adequate states.
\end{theorem}

{\bf Proof:} Every alternating link diagram is adequate, so its
states $s_+$ and $s_-$ are adequate. Note that every
non-alternating link diagram can be transformed into some
alternating diagram and its mirror image by crossing changes which
correspond to changes between positive and negative markers.
Hence, two adequate states of a non-alternating diagram can be
obtained by appropriate choice of markers corresponding to
crossing changes transforming the non-alternating diagram to the
alternating one. $\Box $

The first link that has an adequate state other then $s_+$ and
$s_-$ is the knot $4_1$ ($2\,2$) and it is illustrated in Fig. 6.

\begin{figure}[th]
\centerline{\psfig{file=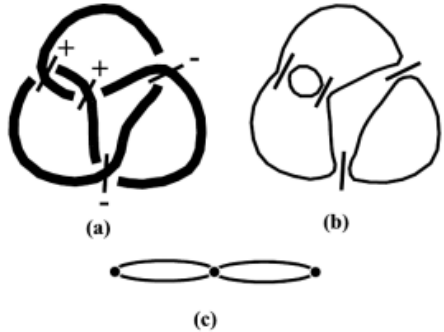,width=2.40in}} \vspace*{8pt}
\caption{(a) Minimal diagram of the figure-eight knot with two
$+$markers, and two $-$markers; (b) state circles; (c) the
associated adequate graph $G_s$. \label{f1.3}}
\end{figure}

The minimal diagram of inadequate knot $2\,0.-3.-2\,0.2$ has as
many as 11 adequate states. First two are obtained from the
alternating diagram $2\,0.3.2\,0.2$ and its mirror image. The
remaining nine adequate states can be obtained from other adequate
diagrams, one corresponding to the minimal diagram
$2\,0.-3.-2\,0.-2$ and the other to the non-minimal diagram
$-2\,0.3.2\,0.-2$ which is reducible to 10-crossing
non-alternating knot $10_{124}$ ($5,3,-2$) (Fig. 7).

\begin{figure}[th]
\centerline{\psfig{file=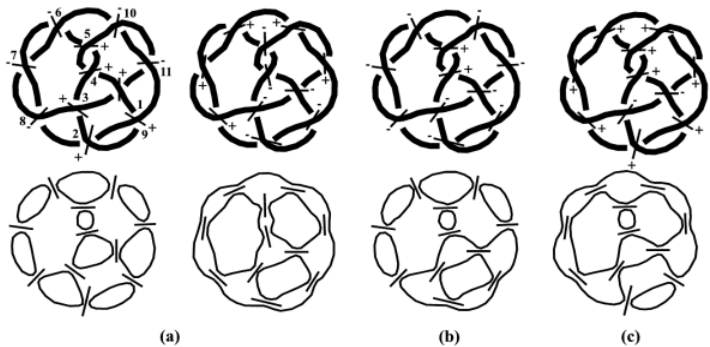,width=3.80in}}
\vspace*{8pt} \caption{(a) Two adequate states of the inadequate
knot diagram $2\,0.-3.-2\,0.2$ obtained from the alternating knot
$2\,0.3.2\,0.2$; (b) adequate state of the same diagram
corresponding to the minimal diagram $2\,0.-3.-2\,0.-2$; (b) its
adequate state corresponding to the non-minimal diagram
$-2\,0.3.2\,0.-2$, which is reducible to 10-crossing
non-alternating knot $10_{124}$ ($5,3,-2$). \label{f1.3}}
\end{figure}

\begin{theorem}
Vertex connectivity of every adequate graph $G_{s_+}$ or $G_{s_-}$
corresponding to an alternating diagram $D$ is greater then $1$.
Vertex connectivity of every adequate graph $G_{s_+}$ or $G_{s_-}$
corresponding to a non-alternating minimal diagram $D$ is $1$.
\end{theorem}

The same statement is not true for adequate graphs $G_{s}$
obtained from other states. For example, the adequate graph $G_s$
of the minimal non-alternating diagram of the knot $10_{155}$ =
$-3:2:2$ (Fig. 8) has the vertex connectivity 4.

\begin{figure}[th]
\centerline{\psfig{file=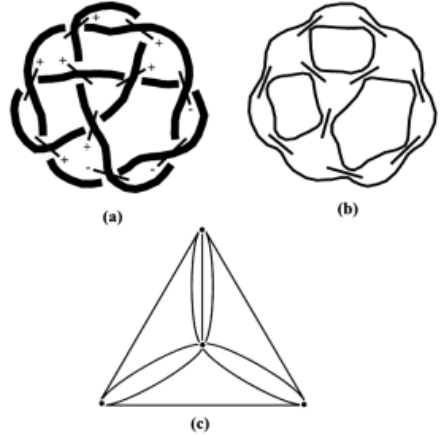,width=2.40in}} \vspace*{8pt}
\caption{(a) Minimal diagram of the knot $10_{155}$ with markers;
(b) state circles; (c) the associated graph $G_s$ with the vertex
connectivity 4. \label{f1.1}}
\end{figure}

\begin{definition}
The minimal number od adequate states taken over all diagrams of a
link $L$ is called the {\it adequacy number} of link $L$ and
denoted by $a(L)$.
\end{definition}

\begin{lemma}
All minimal diagrams of the same alternating link have the same
number of adequate states.
\end{lemma}

Since changing marker in one vertex is equivalent to the crossing
change, we conclude that the number of adequate states is
invariant of a link diagram independent from the signs of
crossings. This means that the number of adequate states is the
same for every alternating diagram and all non-alternating
diagrams obtained from it by crossing changes. Moreover, this can
be generalized to families of links, since adding a bigon to the
chain of bigons does not change the adequacy of a diagram.

\begin{lemma}
The number of adequate states $a(L)$ is the invariant of a family
of alternating links $L$ and it is realized in any minimal diagram
of the link family.
\end{lemma}

\begin{theorem}
The only links with $a(L)=2$ are links of the family $n$
$(n=2,3,4,5,...)$, i.e., $2_1^2$, $3_1$, $4_1^2$, $5_1$,...
\end{theorem}

The numbers of adequate states of two minimal diagrams of a
non-alternating link can be different. The minimal diagram of the
knot $3,2\,1,-2$ has 6 adequate states, and its another minimal
diagram $.2.-2\,0.-1:.-1$ has 8 adequate states, since the source
link of the first diagram is $2,2\,1,2$, and the source link of
the other $.2.2\,0$.

Adequacy numbers of alternating link families obtained from source
links with at most $n=9$ crossings are given in the following
table, where every family is represented by its source link.

\tiny

\bigskip
\noindent \begin{tabular}{|c|c|c|c|c|c|c|} \hline
$n=2$  & $2$ & $$ & $$ & $$ & $$ & $$    \\

$$  & $2$ & $$ & $$ & $$ & $$ & $$    \\ \hline \hline

$n=4$  & $2\,2$ & $$ & $$ & $$ & $$ & $$     \\

$$  & $3$ & $$ & $$ & $$ & $$ & $$    \\ \hline \hline

$n=5$  & $2\,1\,2$ & $$ & $$ & $$ & $$ & $$     \\

$$  & $4$ & $$ & $$ & $$ & $$ & $$     \\  \hline \hline

$n=6$  & $2\,2\,2$ & $2\,1\,1\,2$ & $2,2,2$ & $$ & $$ & $$     \\

$$  & $5$ & $5$ & $5$ & $$ & $$ & $$     \\  \hline \hline

$n=7$  & $2\,1\,2\,2$ & $2\,1\,1\,1\,2$ & $2,2,2+$ & $2\,1,2,2$ & $.2$ &        \\

$$  & $6$ & $7$ & $8$ & $6$ & $7$ & $$     \\ \hline \hline

$n=8$  & $2\,2\,2\,2$ & $2\,1\,2\,1\,2$ & $2\,2\,1\,1\,2$ & $2\,1\,1\,1\,1\,2$ & $2,2,2,2$  & $2\,2,2,2$       \\

$$  & $8$ & $8$ & $8$ & $9$ & $12$ & $7$       \\ \hline

 $$ & $2\,1\,1,2,2$  & $2\,1,2\,1,2$ & $2,2,2++$ & $2\,1,2,2+$ & $(2,2)\,(2,2)$ &  $.2 1$   \\

 $$ & $9$  & $8$ & $9$ & $9$ & $8$  &  $10$   \\  \hline

$$ &   $.2:2$ & $.2.2$ & $.2:2\,0$ & $.2.2\,0$ & $$ &    \\

$$ &  $9$ &  $8$ & $8$ & $8$ & $$ &  \\  \hline \hline

$n=9$  & $2\,2\,1\,2\,2$ & $2\,2\,2\,1\,2$ & $2\,1\,2\,1\,1\,2$ & $2\,2\,1\,1\,1\,2$ & $2\,1\,1\,1\,1\,1\,2$ & $2\,1,2\,1,2\,1$     \\

$$  & $9$ & $10$ & $10$ & $11$ & $12$ & $12$     \\  \hline

$$ &  $2\,1\,2,2,2$  & $2\,2\,1,2,2$ & $2\,1\,1\,1,2,2$ & $2\,1,2,2,2$ & $2\,2,2\,1,2$ & $2\,1\,1,2\,1,2$     \\

$$ &  $10$ & $11$ & $13$ & $9$ & $11$ & $10$    \\ \hline

 $$ & $2\,1,2,2++$ & $2,2,2,2+$ & $2\,2,2,2+1$ & $2\,1\,1,2,2+$ & $2\,1,2\,1,2+$ & $(2\,1,2)\,(2,2)$  \\

$$ &  $10$ & $16$  & $12$ & $13$ & $11$ & $11$  \\ \hline

$$ &  $(2,2)\,(2,2)$ & $(2,2)\,1\,(2,2)$  & $.2\,2$ & $.2\,1\,1$ & $.2\,1:2$ & $.2\,1:2\,0$  \\
$$ &  $10$ & $12$  & $11$ & $13$ & $12$ & $12$  \\ \hline

$$ &  $.2\,1.2\,0$ & $.2.2\,0.2$  & $2:2\,0:2\,0$ & $2\,0:2\,0:2\,0$ & $.2.2.2$ & $2:2:2$  \\
$$ &  $11$ & $10$  & $10$ & $9$ & $10$ & $9$  \\ \hline

$$ &  $.2.2.2\,0$ & $2:2:2\,0$  & $.(2,2)$ & $8^*2$ & $8^*2\,0$ & $$  \\
$$ &  $9$ & $9$  & $14$ & $12$ & $13$ & $$  \\ \hline

\end{tabular}

\normalsize

\bigskip

\section{Adequacy polynomial as an invariant of alternating link
families}

Adequate state graphs corresponding to link diagrams can be used
for creating a polynomial invariant of alternating link families.

\begin{definition}
A {\it cut-vertex} (or articulation vertex) of a connected graph
is a vertex whose removal disconnects the graph [Char]. In
general, a cut-vertex is a vertex of a graph whose removal
increases the number of components [Har]. A graph with no
cut-vertices is called a {\it biconnected graph} [Ski]. A {\it
block} is a maximal biconnected subgraph of a given graph.
\end{definition}

The following transformations will be applied to the adequate state
graphs, till the graph cannot be reduced to a graph with lower
number of vertices:

\begin{itemize}
\item (multiple edge reduction) replace every edge of the multiplicity greater than 2 by a
single edge;\footnote{Since chromatic polynomial of a graph and
graph homology does not recognize mutiple edges, this step is not
necessary for further computations [PrPaSa].}
\item (edge chain collapse) replace maximal part of every chain consisting from edges with
vertices of degree 2 by a new edge connecting the beginning vertex
of the first and ending vertex of the last edge;
\item (block move) every block can be moved along the edges of the remaining part of the
graph.
\end{itemize}

\noindent From every adequate state graph $G$ we obtain the
reduced adequate state graph $\overline G$.

\begin{theorem}
Block move preserves graph torsion and chromatic polynomial of a
graph [PrPaSa].
\end{theorem}

Fig. 9 illustrates reduction of the graph with 16 vertices (Fig.
9a) to the graph with 13 vertices (Fig. 9b), or to its equivalent
graph (Fig. 9c) obtained from it by block moves, which has the
same torsion and chromatic polynomial as the graph (Fig. 9b).

\begin{figure}[th]
\centerline{\psfig{file=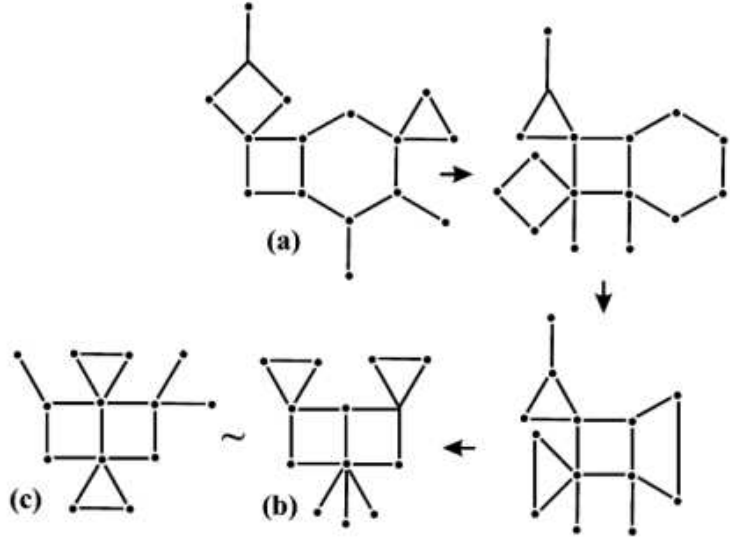,width=4.00in}} \vspace*{8pt}
\caption{Reduction of the graph (a) to the graph (b), and graph
(c) equivalent to (b). \label{f1.1}}
\end{figure}

Consider an arbitrary minimal diagram of an alternating link $L$.
Let $G_i$ denote the corresponding state graphs for all adequate
states of a diagram $D_L$ and $\overline G_i$ reduced state graphs
($i=1,2,...,a(L)$), where $a(L)$ is the adequacy number of $L$ (Def.
8).

\begin{definition}
The {\it adequacy polynomial} of any alternating diagram $D_L$ is a
polynomial in two variables determined by
$A(x,y)=\sum_{i=1}^{a(L)}x^{\overline t_i}\overline P_i(y)$ where
 $\overline P_i(y)=P(\overline G_i)$ denotes the chromatic polynomial of
 a reduced state graph $G_i$ and $\overline t_i$ is the power of
 $\mathbb{Z}$$_2$ torsion of the first chromatic graph
  cohomology $H_{A_m}^{1,h}(G_i)$ over algebra of truncated
  polynomials $A_m=$$\mathbb{Z}$$[x]/x^2=0$ in the grading
  $h=(m-1)(v-2)+1$ where $v$ denotes number of vertices of the graph $G_i$.
\end{definition}

\begin{theorem}
Adequacy polynomial is the same for all minimal diagrams of all
alternating links belonging to the same family, which satisfy the
condition $|a|+k_a\ge 3$\footnote{Please compare this additional
condition with the definition of a family of link diagrams (Def. 5):
according to the additional condition all chains of bigons must be
of the length greater then 2.}.
\end{theorem}

The computation of adequacy polynomial is illustrated on the example
of link $3\,1\,5\,4$. This link has 3 different minimal diagrams:
$3\,1\,5\,4$\footnote{This diagram can be also written as
$(((3,1),1,1,1,1,1),1,1,1,1)$.}, $((1,(1,3),1,1,1,1),1,1,1,1)$, and
$((1,1,(3,1),1,1,1),1,1,1,1)$ (Fig. 10). For the reduced adequate
state graphs $\overline G_i$ ($i=1,2,...,6$) corresponding to the
first minimal diagram, the sequence $(1,2,2,1,1,0)$ represents
powers of $\mathbb{Z}$$_2$-torsion $\overline t_i$ for $m=3$, and
the following is the list of chromatic polynomials:

\medskip

1) $6y-15y^2+14y^3-6y^4+y^5$,

2) $4y-12y^2+13y^3-6y^4+y^5$,

3) $-4y+16y^2-25y^3+19y^4-7y^5+y^6$,

4) $-18y+81y^2-156y^3+168y^4-110y^5+44y^6-10y^7+y^8$,

5) $-2y+5y^2-4y^3+y^4$,

6) $-9y+27y^2- 33y^3+21y^4-7y^5+y^6,$

\medskip

\noindent so the adequacy polynomial is

\small

$$A(3\,1\,5\,4)=
-9y-14xy+27y^2+71xy^2+4x^2y^2-33y^3-146xy^3-12x^2y^3+21y^4+$$
$$163xy^4+13x^2y^4-7y^5-109xy^5-6x^2y^5+y^6+44xy^6+x^2y^6-10xy^7+xy^8.$$

\normalsize

\noindent This polynomial is invariant of link family $p\,1\,q\,r$
($p,q,r\ge 3$).

If we compute the adequacy polynomial from the second or third
diagram, we obtain the same sequence $\overline t_1,\overline
t_2,...,\overline t_6=(1,2,2,1,1,0)$ and the same list of chromatic
polynomials, so the final result remains the same.

All minimal diagrams of the link family $p\,1\,q\,r$ ($p,q,r\ge
3$) have the same adequacy polynomial.

{\bf Conjecture 1.} Adequacy polynomial distinguishes all
alternating link families (up to mutation).

This conjecture is verified for all alternating links with at most
$n=12$ crossings. If the conjecture does not hold in general, one
may consider various adequacy polynomials obtained by taking into
consideration other gradings in first homology or the whole groups
(possibly higher in homology) or changing algebra. Moreover,
depending on the algebra one may consider torsions other then
$\mathbb{Z}$$_2$, if they exist.

Adequacy polynomial of any family of alternating links can be
computed from a minimal diagram of the link $L$ representing this
family with all chains of bigons of the length 3. For subfamilies
we use links with some parameters equal 2, and all the other equal
3. For the general Conway symbol $p\,1\,q\,r$ ($p,q,r\ge 2$), we
need to distinguish the following cases:

\begin{enumerate}
\item $2\,1\,2\,2$ with $A(x,y)=y+2xy-2y^2-5xy^2-4x^2y^2+y^3+6xy^3+8x^2y^3-4xy^4-5x^2y^4+xy^5 +
x^2y^5$;
\item $p\,1\,2\,2$ with $A(x,y)=-8y-8xy+25y^2+26xy^2-4x^2y^2-32y^3-33xy^3+8x^2y^3+21y^4+21xy^4 -
5x^2y^4-7y^5-7xy^5+x^2y^5+y^6+xy^6 $, $p\ge 3$;
\item $2\,1\,q\,2$ with $A(x,y)=-2xy-4x^2y+9xy^2+8x^2y^2-14xy^3-5x^2y^3+11xy^4+x^2y^4-5xy^5+x y^6$, $q\ge 3$;
\item $2\,1\,2\,r$ with $A(x,y)=10xy-8x^3y-27xy^2+28x^3y^2+29xy^3-38x^3y^3-15xy^4+25x^3y^4+3xy^5-
8x^3y^5+x^3y^6$, $r\ge 3$;
\item $p\,1\,q\,2$ with $A(x,y)=-9y^2-2xy^2-4x^2y^2+27y^3+5xy^3+8x^2y^3-33y^4-4xy^4-5x^2y^4+21y^5+
xy^5+x^2y^5-7y^6+y^7$, $p,q\ge 3$;
\item $p\,1\,2\,r$ with $A(x,y)=-9y+6xy+16x^2y+27y^2-17xy^2-60x^2y^2-33y^3+19xy^3+92x^2y^3+21y^4-
10xy^4-75x^2y^4-7y^5+2xy^5+35x^2y^5+y^6-9x^2y^6+x^2y^7$, $p,r\ge
3$;
\item $2\,1\,q\,r$ with $A(x,y)=8xy+8x^2y-20xy^2-32x^2y^2+20xy^3+54x^2y^3-10xy^4-50x^2y^4+2xy^5+
27x^2y^5-8x^2y^6+x^2y^7$, $q,r\ge 3$;
\item $p\,1\,q\,r$ with $A(x,y)=-9y-14xy+27y^2+71xy^2+4x^2y^2-33y^3-146xy^3-12x^2y^3+21y^4+163xy^4
+13x^2y^4-7y^5-109xy^5-6x^2y^5+y^6+44xy^6+x^2y^6-10xy^7+xy^8,$
$p,q,r\ge 3$.
\end{enumerate}

\begin{figure}[th]
\centerline{\psfig{file=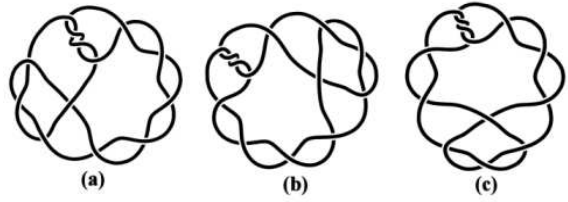,width=3.00in}} \vspace*{8pt}
\caption{Minimal diagrams (a) $3\,1\,5\,4$; (b)
$((1,(1,3),1,1,1,1),1,1,1,1)$; (c) $((1,1,(3,1),1,1,1),1,1,1,1)$.
\label{f1.1}}

\end{figure}

Without any changes, computation of adequacy polynomial can be
extended to families of virtual links. Equivalents of Theorem 21 and
Conjecture 1 hold for alternating virtual links.

The equivalent of Conjecture 1 is verified by computer
calculations for all families of virtual knots derived from real
knots with at most $n=8$ crossings.

\begin{figure}[th]
\centerline{\psfig{file=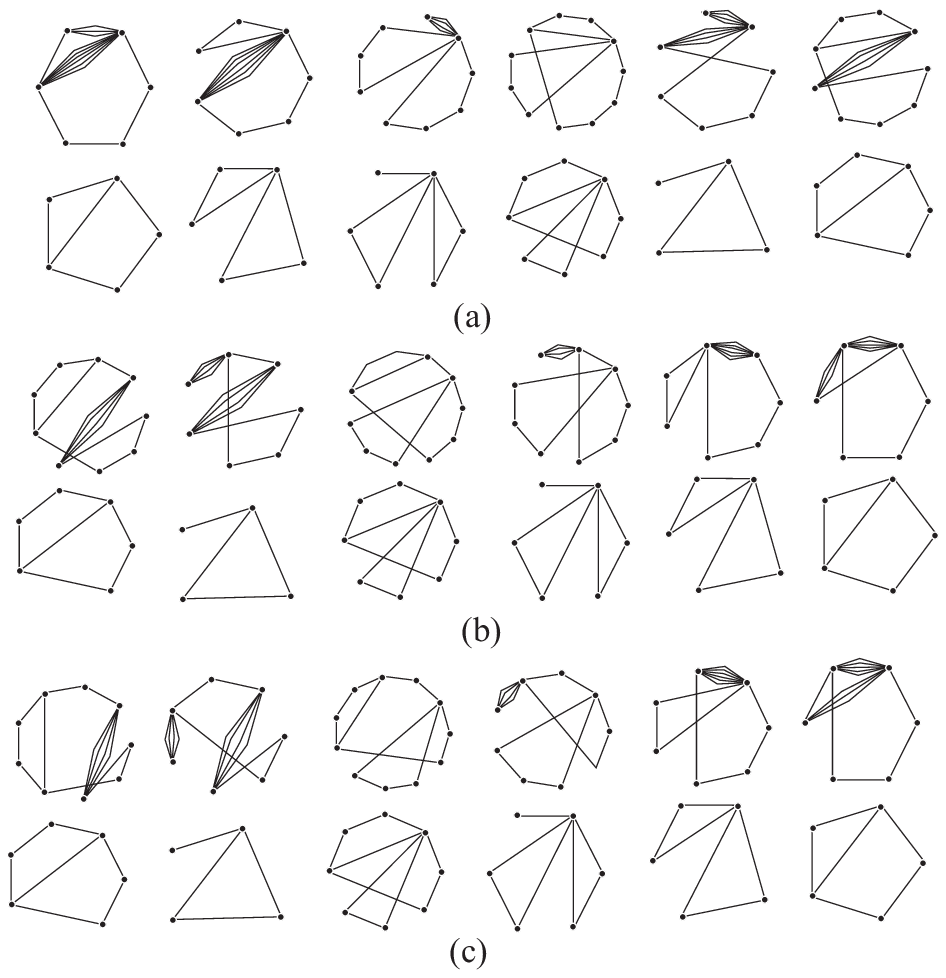,width=4.20in}} \vspace*{8pt}
\caption{Adequate state graphs of the diagrams (a) $3\,1\,5\,4$,
(b) $((1,(1,3),1,1,1,1),1,1,1,1)$, (c)
$((1,1,(3,1),1,1,1),1,1,1,1)$ and their corresponding reduced
graphs. \label{f1.1}}
\end{figure}

Definition of the adequacy polynomial (Def. 10) contains first
chromatic graph homology in the specific grading coming from the
interpretation of Hochschild homology as the chromatic graph
homology of a polygon, i.e., $H_{A_m}^{1,(m-1)(v-2)+1}(G)$ where $G$
is a graph and $v$ denotes the number of its vertices and
$A_m=$$\mathbb{Z}$$[x]/x^m$ for $m \geq 3$. The reason why we have
excluded algebra $A_2$ is that it does not distinguish some
generating links (e.g., $2\,2\,1\,1\,1\,2$ from $2\,1\,1,2\,1,2$).
According to the computations for all generating links with $n\le
12$ crossings, for $3\leq m \leq 5$ adequacy polynomial
distinguishes all families of alternating links with at most $n=12$
crossings (up to mutation).

Notice that adequacy polynomials of the family $3\,1\,3\,3$
computed for $m=2,3,...,8$ are the same, but this is not true on
general:  according to the computer calculations for $2\le m\le
8$, the family $.p$ ($p>2$) will have two different polynomials

$$2y-10xy-10x^2y-4y^2+21xy^2+27x^2y^2+2y^3-14xy^3-31x^2y^3+3xy^4+20x^2y^4-7x^2y^5+x^2y^6$$

\noindent for odd $m$, and

$$2y-4xy-10x^2y-6x^4y-4y^2+10xy^2+27x^2y^2+11x^4y^2+2y^3-8xy^3-31x^2y^3-$$
$$6x^4y^3+2xy^4+20x^2y^4+x^4y^4-7x^2y^5+x^2y^6$$ for even $m$.

In the computation of adequacy polynomials we can also use second
graph homology $H_{A_m}^{2,(m-1)(v-2)}(G)$, which for $m=2,4,6$
distinguishes all alternating link families corresponding to links
with at most $n=12$ crossings.

\bigskip

\noindent {\bf Acknowledgements}

\bigskip

\noindent We would like to express our gratitude to Alexander
Stoimenow and J\' ozef Przyticki for corrections, advice and
suggestions.

\bigskip

\bigskip

\noindent {\bf References}

\bigskip

\small

\noindent [Char] Chartrand, G. (1985) {\it Introductory Graph
Theory}. Dover, New York.

\medskip

\noindent [Cro] Cromwell, P. (2004) {\it Knots and Links}.
Cambridge University Press, Cambridge.

\medskip

\noindent [JaSaz] Jablan, S.~V., Sazdanovi\' c, R. (2007) {\it
LinKnot- Knot Theory by Computer}. World Scientific, New Jersey,
London, Singapore.

\medskip

\noindent [KidSto] Kidwell, M.~E., Stoimenow, A. (2003) Examples
Relating to the Crossing Number, Writhe, and Maximal Bridge Length
of Knot Diagrams. Michigan Math. J. {\bf 51}, 3--12.

\medskip

\noindent [Li] Lickorish, W.B.R. (1991) An Introduction to Knot
Theory. Springer-Verlag, New York, Berlin, Heidelberg.

\medskip

\noindent [LiThi] Lickorish, W.B.R., Thistlethwaite, M. (1988) Some
links with non-trivial polynomials and their crossing-numbers.
Comment. Math. Helvetici {\bf 63}, 527--539.

\medskip

\noindent [PrAs] Przytycki, J.~H., Asaeda, M.~M. (2004) Khovanov
homology: torsion and thickness. arXiv:math/0402402v2 [math.GT]

\medskip

\noindent [PrPaSa] Przytycki, J.~H., Pabiniak, M.~D., Sazdanovic,
R. (2006) On the first group of the chromatic cohomology of
graphs. arXiv:math/0607326v1 [math.GT]

\medskip

\noindent [Pr] Przytycki, J.~H. (2008) Private Communication.

\medskip

\noindent [Ski] Skiena, S. (1990) {\it Implementing Discrete
Mathematics: Combinatorics and Graph Theory with Mathematica}.
Addison-Wesley, Reading, MA.

\medskip

\noindent [Stoi] Stoimenow, A. (2007) Non-triviality of the Jones
polynomial and the crossing numbers of amphicheiral knots.
arXiv:math/0606255v2 [math.GT]

\medskip

\noindent [Stoi1] Stoimenow, A. (2008) Tait's Conjectures and Odd
Crossing Number Amphicheiral Knots. Bulletin (New Series) of the
American Mathematical Society, S 0273-0979(08)01196-8.

\medskip

\noindent [Stoi2] Stoimenow, A. (2008) Private communication.

\medskip

\noindent [Stoi3] Stoimenow, A. (2008) On the crossing number of
semiadequate links.

\medskip

\noindent [Thi] Thistlethwaite, M. (1988) Kauffman polynomial and
adequate links. Invent. Math. {\bf 93}, 2, 258--296.

\normalsize

\end{document}